%% file: IntQForms_ExtRootSubs_II.tex
\documentclass[a4paper,11pt,reqno]{amsart}

\usepackage{graphicx}
\usepackage{amsmath}
\usepackage{amssymb}
\usepackage{amsfonts}
\usepackage{epsfig}
\usepackage{epigraph}
\usepackage[hang]{caption}
\usepackage{hhline}
\usepackage{comment}
\usepackage{mathrsfs}

\usepackage{graphics}
\usepackage{float}

\usepackage{mathtools}

\usepackage[T1]{fontenc}

\graphicspath{{eps_Linkages/}}

 \DeclarePairedDelimiter\abs{\lvert}{\rvert}%

\makeatletter

\newtheorem{theorem}{Theorem}[section]
\newtheorem*{theorem-non}{Theorem}

\@addtoreset{equation}{section}

\def\Sp{{\rm Span}}

\def\rank{{\rm rank}}

\def\botG{\stackrel{}{\Gamma}}

\newtheorem{remark}[theorem]{Remark}
\newtheorem{proposition}[theorem]{Proposition}
\newtheorem*{proposition-non}{Proposition}
\newtheorem{lemma}[theorem]{Lemma}

\def\PerfProof{{\it Proof.\ }}

\begin{document}

\title{Integer quadratic forms and extensions of subsets of
 linearly independent roots}
         \author{Rafael Stekolshchik}

\date{}

\begin{abstract}
  We consider subsets of linearly independent roots in a certain root system $\varPhi$.
  Let $S'$ be such a subset, and let $S'$ be associated with any Carter diagram $\Gamma'$.
  The main question of the paper: what root $\gamma \in \varPhi$ can be added to $S'$
  so that $S' \cup \gamma$ is also a subset of linearly independent roots?
  This extra root $\gamma$ is called the {\it linkage root}.
  The vector $\gamma^{\nabla}$ of inner products  $\{(\gamma,\tau'_i)\mid \tau'_i \in S'\}$ is called the {\it linkage label vector}.
  Let $B_{\Gamma'}$ be the Cartan matrix associated with $\Gamma'$.
  It is shown that $\gamma$ is a linkage root if and only if $\mathscr{B}^{\vee}_{\Gamma'}(\gamma^{\nabla}) < 2$, where $\mathscr{B}^{\vee}_{\Gamma'}$ is a quadratic form with the matrix inverse to  $B_{\Gamma'}$.
  The set of all linkage roots for $\Gamma'$ is called a {\it linkage system} and is denoted by
  $\mathscr{L}(\Gamma')$. The Cartan matrix associated with any Carter diagram $\Gamma'$ is conjugate to the Cartan matrix associated with some Dynkin diagram $\Gamma$ \cite{St23}.
  The sizes of $\mathscr{L}(\Gamma')$ and $\mathscr{L}(\Gamma)$ are the same.
  Let $W^{\vee}$ be the Weyl group of the quadratic form $\mathscr{B}^{\vee}_{\Gamma'}$.
  This group acts on the linkage system and forms several
  orbits. The sizes and structure of orbits for linkage systems $\mathscr{L}(D_l)$ and $\mathscr{L}(D_l(a_k))$ are presented.
\end{abstract}

\maketitle

\input 1Intro.tex

\input 2MainRes.tex

\input 3QuadrForm.tex

\input 4WeylGr_Dual.tex

\input 5DynkinExt.tex

\input 6DTypeLinkages.tex

\input 7ExamplesDiagr.tex

\begin{appendix}
\input App_Gabrielov.tex

\end{appendix}
\input biblio_ext_rs.tex

\end{document}

%% file: 1Intro.tex
\section{\bf Introduction}

\subsection{Brief summary}

The Carter diagram is a generalization of the Dynkin diagram that allows cycles of even length.
Let $\Gamma$ be an arbitrary Carter diagram. 
A set of linearly independent roots $S = \{\tau_1,\dots,\tau_n\}$ is called $\Gamma$-set,
if the vertices of the diagram $\Gamma$ are in one-to-one correspondence with the roots of $S$
and the inner products $(\tau_i, \tau_j)$ correspond to the edges of $\Gamma$.
Precise definitions of what a Carter diagram and $\Gamma$-set are will be given a little further down, 
see $\S\ref{sec_root_syst}$ and $\S\ref{sec_G_set}$.
 
Let $\widetilde\Gamma$ be the diagram obtained from a Carter diagram $\Gamma$ by adding one extra root $\gamma$,
with its bonds, so that the roots corresponding to $\widetilde\Gamma$ form
a linearly independent root subset. We would like to describe such roots in the general case.
The extra root $\gamma$ is called a {\it linkage root},
the diagram $\widetilde\Gamma$ is called a {\it linkage diagram}.
More precise definitions will be given shortly.  

The simple criterion for a given vector $\gamma$ to be a linkage root given in the terms of
the inverse quadratic form $\mathscr{B}^{\vee}$ to the Cartan matrix associated with $\Gamma$ 
is as follows:
  \begin{equation}
     \mathscr{B}^{\vee}_{\Gamma}(\gamma^{\nabla}) < 2,
  \end{equation}  
where $\gamma^{\nabla}$ is so-called {\it linkage label vector},
the vector of the inner products of $\gamma$ with roots from a given $\Gamma$-set
(Theorem \ref{th_B_less_2}).
There exists one-to-one correspondence between 
linkage diagrams and linkage label vectors.     
The Weyl group $W^{\vee}$ preserving the quadratic form $\mathscr{B}^{\vee}$
acts on the space of the linkage diagrams, or, equivalently, 
on the space of linkage label vectors, see \S\ref{sec_Weyl_gr}. 

The set $\mathscr{L}_{\widetilde\Gamma}(\Gamma)$ of all linkage systems 
for the given vertex extension 
$\Gamma \prec \widetilde\Gamma$ is said to be {\it partial linkage system}. 
 (For the definition of vertex extension, see \S\ref{sec_vert_ext}.)
The union of all partial systems using all possible extensions  $\Gamma \prec \widetilde\Gamma$
is called the {\it full linkage system} of $\Gamma$, or, for short, {\it linkage system}:
\begin{equation}
      \mathscr{L}(\Gamma') = \bigcup_{\Gamma ' \prec \widetilde\Gamma} \mathscr{L}_{\widetilde\Gamma}(\Gamma').
\end{equation}
The size and structure of linkage systems $\mathscr{L}(D_l)$  and $\mathscr{L}(D_l(a_k))$ are described in 
Theorem \ref{th_size_Dtype}, 
and Figs.~\ref{Dlpu_linkages_2}--\ref{D6pu_loctets}, Figs.~\ref{D4_loctets}--\ref{D4a1_loctets},
see also Table \ref{tab_number_linkages}.
The group $W^{\vee}$ acts on the linkage system and forms several orbits that are not connected to each other.

The Carter diagrams having the same rank and type form {\it homogeneous classes}. 
For example, $E_6$, $E_6(a_1)$  and $E_6(a_2)$ form the homogeneous class $C(E_6)$.
Any homogeneous class contains only one Dynkin diagram, see \S\ref{sec_homog_class}.
For each pair of diagrams $\{\widetilde\Gamma, \Gamma\}$ out of one homogeneous class,
there exists {\it transition} matrix $M$ mapping each $\widetilde\Gamma$-set $\widetilde{S}$
to some $\Gamma$-set $S$ (Theorem \ref{th_invol}). The matrix $M$ acts only on one element in $\widetilde{S}$,
the remaining elements remain fixed. The transition matrix is an involution. 
Transition matrices make it possible to connect objects associated with different Carter diagrams
from the same homogeneous class, see \S\ref{sec_some_prop}

The transition mappings are close in properties to the transformations  
introduced by A.Gabrielov for systematic study of quadratic forms associated with singularities, 
and to the transformations introduced by S.Ovsienko for the study
weakly positive unit quadratic forms, see \S\ref{sec_other_context}.

\subsection{Root system, Cartan matrix, Carter diagram}
  \label{sec_root_syst}

Let $\varPhi$ be a finite root system, $W$ be the corresponding finite Weyl group.
We consider only root systems with simply-laced Dynkin diagrams.
The Cartan matrix associated with $\varPhi$ is denoted by $B$.
Let $\mathscr{B}$ be the quadratic Tits form associated with $B$ and $(\cdot,\cdot)$
be the inner product induced by $\mathscr{B}$.  For the Cartan matrix $B$,
the following well-known property\footnotemark[1] holds:
\begin{equation}
 \label{eq_Kac}
   \mathscr{B}(\alpha) = 2 \Longleftrightarrow \alpha \in \varPhi.
\end{equation}
\footnotetext[1]{For the  matrix $B/2$, eq. \eqref{eq_Kac} looks like this:
$\mathscr{B}(\alpha) = 1 \Longleftrightarrow \alpha \in \varPhi$, see \cite[Prop. 1.6]{K80}.}
In \cite{Ca72}, Carter introduced admissible diagrams to describe conjugacy classes in $W$.
These diagrams can be used to describe other objects,
such as root subsets\footnotemark[2] and quadratic forms.
Let $S \subset \varPhi$ be a root subset, roots of $S$ are not necessarily simple.
\footnotetext[2]{In what follows, the phrase ``root subset'' means a subset of
 linearly independent roots.}
To the subset $S$ we associate some diagram $\Gamma$ that provides one-to-one correspondence
between roots of $S$ and nodes of $\Gamma$.
The diagram $\Gamma$ is said to be {\it admissible} if the following two conditions hold:
\begin{equation}
 \label{eq_def_adm}
  \begin{split}
  & (a)  \text{ The nodes of } \Gamma \text{ correspond to a subset of linearly independent roots in } \varPhi. \\
  & (b) \text{ If a subdiagram of } \Gamma \text{ is a cycle, then it contains an even number of nodes.}
  \end{split}
\end{equation}
 In admissible diagrams all edges are solid, whereas Carter diagrams have solid and
 dotted diagrams that distinguish negative and positive inner products, see \S\ref{sec_adm_diagr}.
 In \cite{St17}, it was observed that the cycles in the Carter diagrams
 contain at least one solid edge $\{\alpha_1, \beta_1\}$ (on which $(\alpha_1, \beta_1) < 0$)
 and at least one dotted edge $\{\alpha_2, \beta_2\}~$(on which $(\alpha_2, \beta_2) > 0$),
 see \S\ref{sec_solid_dotted}. In this paper we introduce new objects - linkage roots, linkage diagrams and others,
 which are also closely related to Carter diagrams.

\subsection{The partial Cartan matrix}
 \label{sec_G_set}
Let $\Gamma$ be a Carter diagram.
A set of linearly independent roots $S = \{\tau_1,\dots,\tau_n\}$ is called $\Gamma$-set,
if the vertices of the diagram $\Gamma$ are on one-to-one correspondence with the roots of $S$
and the inner products $(\tau_i, \tau_j)$ correspond to the edges of $\Gamma$ as follows:
 \begin{equation*}
   (\tau_i, \tau_j) =
     \begin{cases}
        0, & \text {if } \tau_i \text{ and } \tau_j \text{ are not connected}, \\
        -1, & \text{if edge } \{\tau_i, \tau_j\} \text{ is solid}, \\
        1, & \text{if edge } \{\tau_i, \tau_j\} \text{ is dotted}. \\
     \end{cases}
 \end{equation*}
 Besides this, $(\tau_i, \tau_i) = 2$ for any $i$.
 Similarly to the Cartan matrix associated with Dynkin diagrams, for each pair \{$\Gamma$, S\},
 we determine the following matrix:

 \begin{equation}
   \label{canon_dec_2}
   B_{\Gamma} :=
      \left (
        \begin{array}{cccccc}
         (\tau_1, \tau_1) & \dots & (\tau_1, \tau_n) \\
         \dots                & \dots & \dots \\
         (\tau_n, \tau_1) & \dots & (\tau_n, \tau_n) \\
        \end{array}
      \right ).
 \end{equation}
 ~\\
 We call the matrix $B_{\Gamma}$ the {\it partial Cartan matrix} corresponding to $\Gamma$.
 Why partial? Since in this case some off-diagonal elements are $1$, they  correspond
 to the dotted edges.
 The partial Cartan matrix $B_{\Gamma}$ is well-defined since products $(\tau_i, \tau_j)$ in \eqref{canon_dec_2}
 do not depend on the choice of the $\Gamma$-set $S$. The elements of the partial
 Cartan matrix are uniquely determined by the diagram $\Gamma$.

 \subsection{Some properties of Dynkin and Carter diagrams}
   \label{sec_some_prop}
\begin{figure}[!ht]
\centering
\includegraphics[scale=0.33]{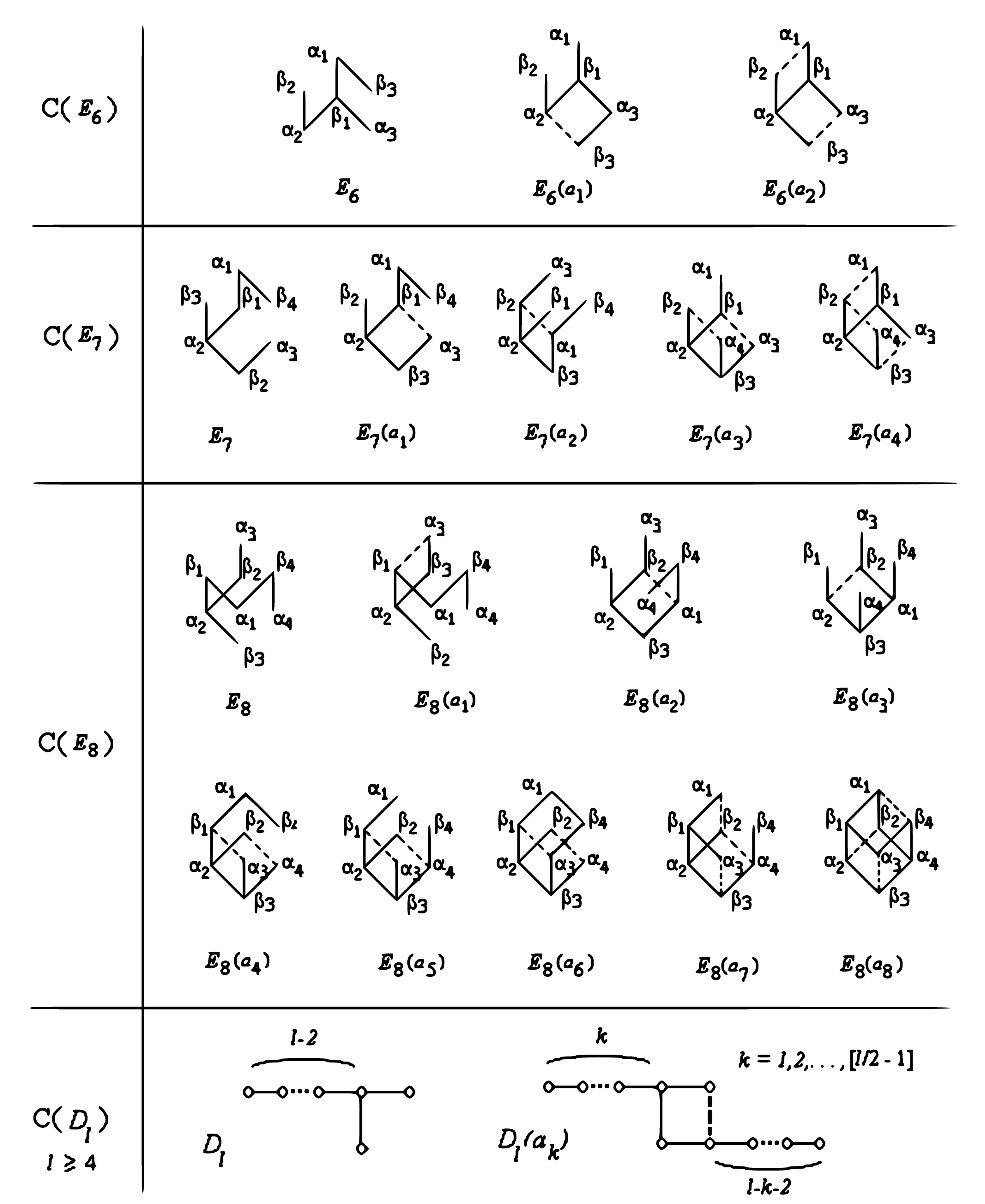}
\caption{Homogeneous classes of Carter diagrams $C(D_l), l \geq 4$ and $C(E_l), l=6,7,8$.}
\label{fig_all_diagr_ED}
\end{figure}

 Any simply-laced Dynkin diagram is the particular case of the Carter diagram. We will use the phrase
 ``Carter diagram'' to refer to a Carter diagram which is not a Dynkin diagram unless otherwise stated.
 For a list of all Carter diagrams, see Fig.~\ref{fig_all_diagr_ED}.
 This list is divided into homogeneous classes, see \S\ref{sec_homog_class}.

  (i) Consider the orbit of any root subset  associated with a Carter diagram (resp.
 Dynkin diagram) under the action of the Weyl group.
 In the case of a Carter diagram, each root subset in this orbit contains non-simple roots.
 In the case of a Dynkin diagram, some root subset on this orbit consists only of simple roots.

 (ii) For both Carter and Dynkin diagrams, the corresponding quadratic forms
 are positive definite; for both of them the spectrum of the corresponding Cartan matrix lies in the
 open interval $(0, 4)$, see Propositions \ref{prop_restr_forms} and \ref{prop_spectr}.
 The latter follows from the fact that the partial Cartan matrix
 $B_{\Gamma'}$ for the Carter diagram
 $\Gamma'$ is related to the Cartan matrix $B_{\Gamma}$ for some Dynkin diagram $\Gamma$
 by a certain transition matrix $M$:
\begin{equation}
  \label{eq_matr_M}
    {}^t{M} \cdot B_{\Gamma'} \cdot M = B_{\Gamma},
\end{equation}
see \S\ref{sec_trans_mappings}.

 (iii) The transition matrix $M$ modifies the root subset, the quadratic form as \eqref{eq_matr_M}
 and even the underlying diagram.
 The transition matrix combines the Carter diagrams and the corresponding Dynkin diagram
 together into one homogeneous class, see \S\ref{sec_homog_class}. Each homogeneous class
 contains only one Dynkin diagram.

 (iv) The Carter diagrams contain $4$-cycles $D_4(a_1)$ as subdiagrams. In principle, Carter diagrams can contain
  the cycles of length > $4$. It is shown in \cite[Th 3.1]{St17} that any Carter diagram containing
  $l$-cycles, where $l >4$, is equivalent to another Carter diagram containing only $4$-cycles.
  To realize this equivalence, an explicit transformation was constructed  that maps each Carter
  diagram with long cycles into some Carter diagram containing only $4$-cycles.

  (v) The correspondence between root subsets and Carter diagrams is not one-to-one.
  For example, the set of four linear independent roots $\{\alpha_1, \alpha_2, \beta_1,  \beta_2\}$
  connected as shown in Fig.~\ref{D4a1_loctets},
  can be associated with exactly two conjugacy classes having the following representatives:
  $w_o = s_{\alpha_1}s_{\beta_1}s_{\alpha_2}s_{\beta_2}$ (associated with $D_4$) and
  $w_1 = s_{\alpha_1}s_{\alpha_2}s_{\beta_1}s_{\beta_2}$ (associated with $D_4(a_1)$),
  \cite[\S\S1.2.3--1.2.4]{St17}. The element $w_o$ corresponds to the Coxeter transformation class,
  and $w_1$ corresponds to the non-Coxeter class.
  As Vavilov and Migrin noted in \cite[\S1.2]{VM21},
  this single phenomenon explains  all non-Coxeter classes for all
  Carter diagrams with cycles.

\subsection{Homogeneous Carter diagrams and the transition matrix}
  \label{sec_homog_class}
  The Dynkin diagrams $A_l$, where $l \ge 1$ (resp. $D_l$,  $l \ge 4$; resp. $E_l$, $l = 6,7,8$)
  are said to be the {\it Dynkin diagrams of $A$-type}
  (resp. {\it $D$-type}, resp. {\it $E$-type}). The Carter diagrams $D_l$, $D_l(a_k)$,
  $l \ge 4$, $1 \leq k \leq \big [ \tfrac{l-2}{2} \big ]$
  (resp. $E_l$, $E_l(a_k)$, $l = 6,7,8$, $k$ are given in \eqref{eq_one_type})
  are said to be the {\it Carter diagrams of $D$-type}, resp. {\it $E$-type}).

 The Carter diagrams of the same type and the same index constitute a {\it homogeneous class} of Carter diagrams.
 Denote by $C(\Gamma)$ the homogeneous class containing the Carter diagram $\Gamma$,
 see \eqref{eq_one_type} and Fig.~\ref{fig_all_diagr_ED}.
\begin{equation}
 \label{eq_one_type}
  \begin{array}{ll}
      & C(E_6) =  \{  E_6, E_6(a_k), k = 1,2 \}, \\
      & C(E_7) =  \{  E_7, E_7(a_k), 1 \leq k \leq 4 \}, \\
      & C(E_8) =  \{  E_8, E_8(a_k), 1 \leq k \leq 8 \}, \\
      & C(D_l) =  \{  D_l, D_l(a_k), 1 \leq k \leq \big [\tfrac{l-2}{2} \big ] \}, \text{ where } l \geq 4.
  \end{array}
\end{equation}

Let $\widetilde{S}$ (resp. $S$) be a $\widetilde\Gamma$-set (resp. $\Gamma$-set).
In \cite[Tables 2-5]{St23},
the {\it transition matrices} $M : \widetilde{S} \longmapsto S$ are constructed for the following
homogeneous pairs  $\{ \widetilde\Gamma, \Gamma \}$:

  \begin{equation}
    \label{eq_pairs_trans}
     \begin{array}{ll}
          & (1)\; \{D_4(a_1), D_4\} \\
          & (2)\; \{D_l(a_k), D_l\} \\
          & (3)\; \{E_6(a_1), E_6\} \\
          & (4)\; \{E_6(a_2), E_6(a_1)\} \\
          & (5)\;  \{E_7(a_1), E_7\}  \\
          & (6)\;  \{E_7(a_2), E_7\}  \\
          & (7)\;  \{E_7(a_3), E_7(a_1)\} \\
          & (8)\;  \{E_7(a_4), E_7(a_3)\} \\
     \end{array}\qquad \qquad
     \begin{array}{ll}
          & ~(9) \;  \{E_8(a_1), E_8\} \\
          & (10)\;  \{E_8(a_2), E_8\} \\
          & (11)\;  \{E_8(a_3), E_8(a_2)\} \\
          & (12)\;  \{E_8(a_4), E_8(a_1)\} \\
          & (13)\;  \{E_8(a_5), E_8(a_4)\} \\
          & (14)\;  \{E_8(a_6), E_8(a_4)\} \\
          & (15)\;  \{E_8(a_7), E_8(a_5)\} \\
          & (16)\;  \{E_8(a_8), E_8(a_7)\} \\
     \end{array}
  \end{equation}
\smallskip

  To our considerations, let us add {\it diagram similarity mapping} $L_{\tau_i}$
  reflecting the root $\tau_i \in \widetilde{S}$ and fixing the remaining roots in $\widetilde\Gamma$-set:
\begin{equation}
     \label{eq_similarity_0}
     L_{\tau_i}: \tau_i \longmapsto -\tau_i.
\end{equation}
  Eq.~\eqref{eq_similarity_0} means that any solid (resp. dotted) edge with vertex $\tau_i$ is mapped
  to the dotted (resp. solid) edge.
  Two Carter diagrams obtained from each other by a sequence of reflections \eqref{eq_similarity_0},
  are called {\it similar Carter diagrams}. The diagrams obtained as images
  of the transition matrices $M$ in \eqref{eq_matr_M} are considered up to similarity of Carter diagrams,
  see \cite[\S 1.3]{St17}.
  The list \eqref{eq_pairs_trans} is called the {\it adjacency list}.

\subsection{Transitions}
  \label{sec_trans_mappings}
The transition matrices proposed in \cite{St23} transform root subsets
associated with the corresponding Carter diagrams. Given a Carter diagram $\Gamma$,
each transition  changes only one root in the subset associated with $\Gamma$.

\begin{theorem}[Transitions theorem] \cite[Th. 4.1]{St23}
 \label{th_invol}
For each pair of diagrams $\{\widetilde\Gamma, \Gamma \}$ out of the list \eqref{eq_pairs_trans},
there exists the \underline{transition matrix} $M$ mapping each $\widetilde\Gamma$-set $\widetilde{S}$  to
some $\Gamma$-set $S$.
The matrix $M$ acts only on one element $\widetilde\alpha \in \widetilde{S}$
and does not change remaining elements in $\widetilde{S}$:
\begin{equation}
 \label{th_transition}
  \begin{cases}
      M \tau_i = \tau_i \text{ for all }
           \tau_i \in \widetilde{S}, \tau_i \neq \widetilde\alpha, \\
      M \widetilde\alpha = \alpha =
         -\widetilde\alpha + \sum t_i \tau_i,
         \text{ the sum is taken over all }
           \tau_i \in \widetilde{S} \text{ except } \widetilde\alpha; ~t_i \in \mathbb{Z}. \\
  \end{cases}
\end{equation}
  The image $S = M\widetilde{S}$  the $\Gamma$-set.
  The transition matrix $M: \widetilde{S} \longmapsto S$ is an
  \underline{involution}:
\begin{equation*}
   M \widetilde\alpha = \alpha \text{ and } M\alpha = \widetilde\alpha.
\end{equation*}
\end{theorem}

\medskip
The transition matrix change the underlying diagram for given subset of roots.
For more details about the transitions matrices $M$ corresponding to the adjacency list \eqref{eq_pairs_trans},
see Tables 2-5 in  \cite[\S4.6]{St23}.
In \S\ref{sect_Cartan} we use transitions $M$ to relate partial Cartan matrices associated with
homogeneous Carter diagrams, as well as to relate linkage systems associated with them, see \S\ref{sec_linkage}.

\subsubsection{Transition mappings in other contexts}
  \label{sec_other_context}
I would like to highlight two areas that use transformations close to the transition matrices from
Theorem \ref{th_invol}.
In 1973, A.~Gabrielov introduced transformations
for the systematic study of quadratic forms associated with singularities, \cite{G73},
see Gabrielov's examples in \S\ref{sec_Gabrielov}.
In 1978, S.~Ovsienko introduced so-called inflation technique for study
the weakly positive unit quadratic forms, \cite{O78}. For the Ovsienko theorem, see \S\ref{sec_Ovsienko}.

\subsection{Admissible diagrams and Carter diagrams}
   \label{sec_adm_diagr}

Each element $w \in W$ can be expressed in the form
 \begin{equation}
   \label{any_roots_0}
    w  = s_{\alpha_1} s_{\alpha_2} \dots s_{\alpha_k}, \text{ where } \alpha_i \in \varPhi \text{ for all } i.
 \end{equation}

Carter proved that $k$ in the decomposition \eqref{any_roots_0} is the smallest if
and only if the subset of roots $\{\alpha_1, \alpha_2, \dots , \alpha_k\}$ is linearly independent;
such a decomposition is said to be {\it reduced}. The admissible diagram
corresponding to the given element $w$ is not unique, since the reduced
decomposition of the element $w$ is not unique.

 Denote by $l_C(w)$ the smallest value $k$
 corresponding to any reduced decomposition \eqref{any_roots_0}.
 The corresponding set of roots $\{\alpha_1,\alpha_2,\dots,\alpha_k\}$ consists of linearly
 independent and not necessarily simple roots, see Lemma \ref{lem_lin_indep}. 
 If $l(w)$ is the smallest value $k$ in \eqref{any_roots_0} such that all roots $\alpha_i$ are {\it simple},
 then $l_C(w) \leq l(w)$.

 \begin{lemma}{\rm{\cite[Lemma 3]{Ca72}}}
  \label{lem_lin_indep}
   Let $\alpha_1, \alpha_2, \dots, \alpha_k \in \varPhi$.
   Then, $s_{\alpha_1} s_{\alpha_2} \dots s_{\alpha_k}$ is reduced
   if and only if $\alpha_1, \alpha_2, \dots, \alpha_k$ are linearly
   independent. \qed
 \end{lemma}

The admissible diagram may contain cycles, since the roots of $S$ are non necessarily simple,
\cite{St17}. Let us fix some basis of roots corresponding to the given
admissible diagram $\Gamma$:
\begin{equation}
\label{eq_alpha_bet}
  S = \{\alpha_1,\dots,\alpha_k,\beta_1,\dots,\beta_h\}
\end{equation}
\noindent
According to (\ref{eq_def_adm}$(a)$), there exists a certain set \eqref{eq_alpha_bet}
of linearly independent roots.Thanks to (\ref{eq_def_adm}$(b)$),
there exists a partition $S = S_{\alpha} \coprod S_{\beta}$
which is said to be the {\it bicolored partition}.
The admissible diagram is bicolored, i.e., the set of nodes can be partitioned
into two disjoint subsets $S_\alpha = \{\alpha_1,\dots,\alpha_k \}$ and $S_\beta = \{\beta_1,\dots,\beta_h\}$,
where roots of $S_\alpha$ (resp. $S_\beta$) are  mutually orthogonal.

 Let $w = w_1 w_2$ be the decomposition of $w$ into the product of two involutions.
 By \cite[Lemma~5]{Ca72} each of $w_1$ and $w_2$ can
 be expressed as a product of reflections as follows:
 \begin{equation}
   \label{two_invol}
      w = w_1{w}_2, \quad \text{ where } \quad
      w_1 = s_{\alpha_1} s_{\alpha_2} \dots s_{\alpha_k}, \quad
      w_2 = s_{\beta_1} s_{\beta_2} \dots s_{\beta_h},
 \end{equation}
 where subset $S_{\alpha} = \{\alpha_1,\dots,\alpha_k\}$ (resp. $S_{\beta} = \{\beta_1,\dots,\beta_h\}$)
 consists of mutually orthogonal roots.
 The decomposition \eqref{two_invol} is said to be a {\it bicolored decomposition}.
 The subset of roots corresponding to $w_1$ (resp. $w_2$) is said to be {\it $\alpha$-set} (resp. {\it $\beta$-set}):
 \index{bicolored decomposition}
\begin{equation}
   \label{two_sets}
       \alpha\text{-set} = \{ \alpha_1, \alpha_2, \dots, \alpha_k \}, \quad
       \beta\text{-set} =  \{ \beta_1, \beta_2, \dots, \beta_h \}.
 \end{equation}

 The element $w$ in \eqref{two_invol} is said to be {\it semi-Coxeter element};
 it represents the conjugacy class associated with
 the admissible diagram $\Gamma$  and root subset $S$.

  Any admissible diagram $\Gamma$ is said to be a {\it Carter diagram}
  if any edge connecting a pair of roots $\{\alpha, \beta\}$ with
  inner product $(\alpha, \beta) > 0$ (resp. $(\alpha, \beta) < 0$)
  is drawn as dotted (resp. solid) edge. There is no edge for the inner product $(\alpha, \beta) = 0$.

 Up to solid/dotted edges, the classification of Carter diagrams coincides with the
 classification of admissible diagrams.

\subsection{Linkage root and linkage diagram}
  \label{sec_linkage}
We consider the {\it extension} of the root subset $S$
from the bicolored partition \eqref{eq_alpha_bet}
by means of the root $\gamma \in \varPhi$,
so that the set of roots
 \begin{equation*}
    S' = \{ \alpha_1, \dots, \alpha_k, \beta_1, \dots, \beta_h, \gamma \}
 \end{equation*}
 is linearly independent. The new diagram $\Gamma'$ is obtained by addition new edges,
 these edges are {\it solid} (resp. {\it dotted})
  for $(\gamma, \tau) = -1$ (resp. $(\gamma, \tau) = 1)$,
  where $\tau \in S$.
~\\
~\\
\begin{minipage}{10.6cm}
 \hspace{5mm}The diagram $\Gamma'$ is said to be
 a {\it linkage diagram} and the root $\gamma$  is said to be a {\it linkage root}.
 The roots $\tau$ corresponding to the new edges ($(\gamma, \tau) \neq 0$) are said to be
 {\it endpoints} of the linkage diagram.
 Consider vector $\gamma^{\nabla}$ defined by \eqref{eq_dual_gamma}.
 This vector is said to be the {\it linkage label vector}.
 There is a one-to-one correspondence between linkage label vectors $\gamma^{\nabla}$
 (with labels $\gamma^{\nabla}_i \in \{-1, 0, 1\}$) and
 linkage diagrams such that $(\gamma, \tau) \in \{-1, 0, 1\}$.
 \end{minipage}
\begin{minipage}{5.4cm}
  \quad
 \begin{equation}
   \label{eq_dual_gamma}
   \gamma^{\nabla} :=
   \left (
    \begin{array}{c}
    (\gamma, \alpha_1) \\
    \dots,     \\
    (\gamma, \alpha_k)  \\
    (\gamma, \beta_1)  \\
    \dots,  \\
    (\gamma, \beta_h) \\
    \end{array}
   \right )
 \end{equation}
\end{minipage}
~\\

\subsection{Solid and dotted edges of linkage diagrams}
 \label{sec_solid_dotted}
 \begin{lemma}\cite[Lem. A.1]{St17}
  \label{lem_edges}
   Every cycle in the Carter diagram contains at least one solid edge and at least one dotted edge.
 \end{lemma}
\begin{figure}[h]
\centering
\includegraphics[scale=0.25]{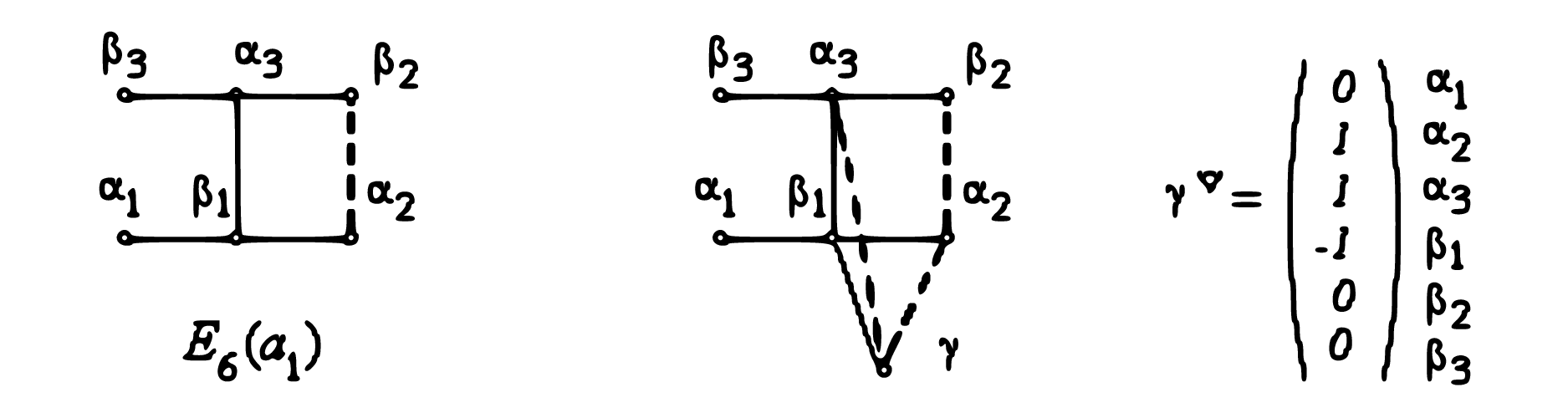}
 \caption{Example of linkage diagram for $E_6(a_1)$.}
\label{fig_E6a1_linkages}
\end{figure}
\noindent
This lemma is also valid for cycles of any linkage diagram.
In Fig.~\ref{fig_E6a1_linkages} if $\{\beta_1, \gamma\}$ is solid then
edge $\{\alpha_3, \gamma\}$ is necessarily dotted, see triangle $\{\alpha_3, \gamma, \beta_1\}$.
Similarly, from triangle $\{\alpha_2, \gamma, \beta_1\}$  the edge $\{\alpha_2, \gamma\}$ is necessarily dotted.

\subsection{Getting rid of triangles}
\begin{figure}[!ht]
\centering
\includegraphics[scale=0.30]{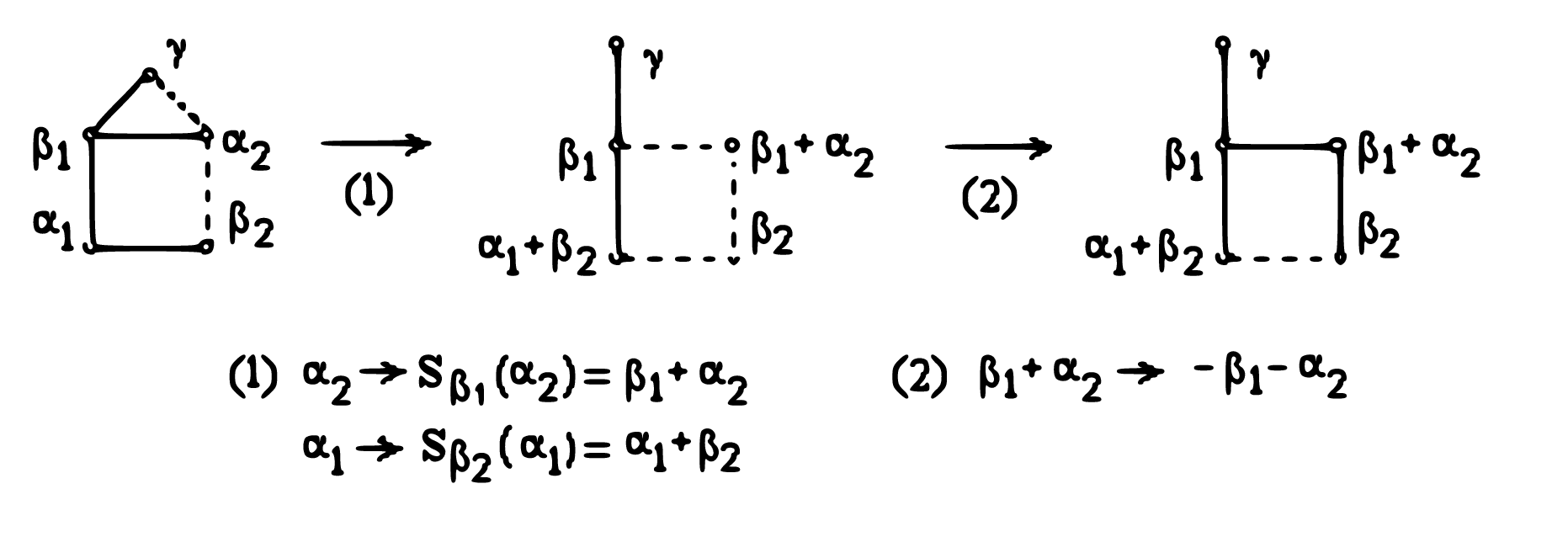}
 \caption{ Eliminating the triangle}
\label{fig_D4a1_triangle}
\end{figure}

By adding the linkage root to a Carter diagram, we get a new diagram
which is not necessarily a Carter diagram. After this, we act by the Weyl group
of the inverse quadratic form and obtain the orbit of the linkage diagrams.
In the most of cases, this orbit contains a Carter diagram,
see Figs.~\ref{D4a1_loctets},~\ref{Dlpu_linkages_2},~\ref{Dk_al_linkages}.
Sometimes the orbit does not contain a Carter diagram because
every linkage diagram on the orbit contains a triangle or other odd cycle.
By changing the basis we can eliminate a triangle or other odd cycle.
The procedure for eliminating the triangle is shown
in Fig.~\ref{fig_D4a1_triangle}.

\subsection{Nilpotent orbits in the simple Lie algebras}
There is a connection (which was pointed out to me by A.~Elashivili) between
 the representatives of nilpotent orbits in simple Lie algebras  and root subsets 
of linearly independent roots represented by Carter diagrams.
By theorem of Jacobson-Morozov, this connection is also a connection between
simple three-dimensional $S$-subalgebras  (or $sl2$-triples) listed by Dynkin,
and root subsets of linearly independent roots, see 
\cite{Elk72}, \cite{Gr08}, \cite{D52}, \cite{On00}, \cite{E75}. \cite{MK21} .

\subsection*{Acknowledgements} I thank A.Elashvili and M.Jibladze for their attention to this paper
and useful discussions during its preparation.

%% file: 2MainRes.tex
\section{\bf The main results}

\subsection{When a vector is a linkage root?}
  \label{sec_linkage_diagr}

  According with \S\ref{sec_linkage}, the linkage diagram is obtained from a Carter diagram $\Gamma$
  by adding some linkage root $\gamma$,
  with its bonds, so that the roots corresponding to vertices of $\Gamma$ together with $\gamma$ form
  a linearly independent root subset.
  We would like to describe such roots in the general case, namely:
\begin{equation}
\label{eq_linkage}
 \begin{array}{l}
 \quad
  \textit{ Let } \Gamma \textit{ be  a Carter diagram, and $S$ some $\Gamma$-set}.
  \textit{ What root $\gamma$ can be added to $S$ } \\
  \textit{ so that $S \cup \gamma$ is a set of linearly independent roots in some
  finite root system? }
  \end{array}
\end{equation}
  With every linkage diagram we associate the linkage label vector.
  There is one-to-one correspondence between linkage label vectors and linkage diagrams, see \S\ref{sec_linkage}.
  The linkage diagram and their linkage label vector $\gamma^{\nabla}$ for the Carter diagram $E_6(a_1)$
  is depicted in Fig.~\ref{fig_E6a1_linkages}. For more examples of linkage diagrams and linkage label vectors
  see Fig.~\ref{D4a1_loctets} (case $D_4(a_1)$), Figs.~\ref{Dk_al_linkages} (case $D_k(a_l)$).

  It turns out that the answer to question \eqref{eq_linkage}
  is very simple in terms of the quadratic form associated with the partial Cartan matrix.
  Let $\gamma$ be a linkage root for $\Gamma$,
 $\gamma^{\nabla}$ be the linkage label vector, see \eqref{eq_dual_gamma}.
 We denote by $L$ the space spanned by roots of $S$, and by $L(\gamma)$ the space spanned by $L$ and $\gamma$.
 Let us write
\begin{equation}
  \label{L_extended_gamma}
   L =  [\tau_1,\dots,\tau_l], \qquad
   L(\gamma) =  [\tau_1,\dots,\tau_l,\gamma].
\end{equation}
 Note that $L(\gamma)$ corresponds to a certain linkage diagram.
 Let $\mathscr{B}_{\Gamma}$ be the quadratic form associated with the partial Cartan matrix $B_{\Gamma}$.
 Then,
\begin{equation*}
  \mathscr{B}_{\Gamma}(\gamma) = \langle B_{\Gamma}\gamma, \gamma \rangle =
  \langle \gamma^{\nabla}, B_{\Gamma}^{-1}\gamma^{\nabla} \rangle,
\end{equation*}
where $\langle\cdot, \cdot\rangle$ is the usual dot product.
Since $B_{\Gamma}$ is positive definite, the eigenvalues of $B_{\Gamma}$
are positive. Therefore, the matrix $B^{-1}_{\Gamma}$ is also positive definite.
We call the quadratic form $\mathscr{B}^{\vee}_{\Gamma}$ corresponding to the matrix
$B^{-1}_{\Gamma}$  the {\it inverse quadratic form}. The form
 $\mathscr{B}^{\vee}_{\Gamma}$ is positive definite.

\subsection{Criterion that a vector is a linkage root}

\begin{theorem-non}[Theorem \ref{th_B_less_2}, Linkage root criterion]
     {\rm(i)} Let $\theta^{\nabla}$ be the linkage label vector
     corresponding to a certain root $\theta \in L(\gamma)$.
     The root $\theta$ is a linkage root
     if and only if
     \begin{equation}
       \label{eq_p_less_02}
          \mathscr{B}^{\vee}_{\Gamma}(\theta^{\nabla}) < 2.
     \end{equation}

    {\rm(ii)} Let $\theta \in L(\gamma)$ be a root with only
     one endpoint $\tau_i \in S$. Then $\theta$ is a linkage root if and only if
     \begin{equation*}
          b^{\vee}_{i,i} < 2,
     \end{equation*}
     where $b^{\vee}_{i,i}$ is the $i$th diagonal element of $B^{-1}_{\Gamma}$.
\end{theorem-non}

Note that criterion \eqref{eq_p_less_02} does not depend on the choice of the $\Gamma$-set $S$,
but only on the quadratic form $\mathscr{B}$ or, what the same, on the diagram $\Gamma$.

Let $\{\Gamma', \Gamma\}$ be a pair of Carter diagrams from the  list \eqref{eq_pairs_trans},
$M$ be the transition matrix transforming $\Gamma'$-set $S'$ to $\Gamma$-set $S$, see \S\ref{sec_trans_mappings},
and let $B_{\Gamma'}$ and $B_{\Gamma}$ be the partial Cartan matrices for $\Gamma'$ and $\Gamma$.
Then, $M$ relates the partial Cartan matrices $B_{\Gamma'}$ and $B_{\Gamma}$ as follows:
\begin{equation*}
   {}^t{M} \cdot B_{\Gamma'} \cdot M = B_{\Gamma}. \\
\end{equation*}

It is shown in {\bf Proposition \ref{prop_quadr_1}} that vector $\gamma^{\nabla}$
is the linkage label vector for $\Gamma$ if and only if \, ${}^t{M}\gamma^{\nabla}$
is the linkage label vector for $\Gamma'$. 

\subsection{Full linkage systems $\mathscr{L}(D_l)$ and $\mathscr{L}(D_l(a_k))$}

Denote by $\mathscr{L}_{\widetilde\Gamma}(\Gamma', S')$ the set of linkage diagrams
\begin{equation}
  \label{eq_part_linkage_syst}
      \mathscr{L}_{\widetilde\Gamma}(\Gamma', S') =
      \{ \gamma^{\nabla} \mid \gamma \in \varPhi(\widetilde\Gamma), \gamma \not\in L(S') \},
\end{equation}
where $\varPhi(\widetilde\Gamma)$ is the root system associated with $\widetilde\Gamma$,
$S'$ some $\Gamma'$-set and $L(S')$ is the linear span $\Sp(S')$. It was shown in
{\bf Proposition \ref{prop_not_dep_S}} that $\mathscr{L}_{\widetilde\Gamma}(\Gamma', S')$ does not
depend on choosing $\Gamma'$-set $S'$. Then, the set \eqref{eq_part_linkage_syst} can be denoted
by $\mathscr{L}_{\widetilde\Gamma}(\Gamma')$.
The set $\mathscr{L}_{\widetilde\Gamma}(\Gamma')$
is said to be {\it partial linkage system $\widetilde\Gamma$ over $\Gamma'$}.
The pair $\{\Gamma', \widetilde\Gamma\}$, where $\Gamma' \in C(\Gamma)$,
is said to be the {\it vertex extension} and is denoted by $\Gamma ' \prec \widetilde\Gamma$.
The union of all partial systems by all possible vertex extensions
is said to be the {\it linkage system} or {\it full linkage system}  of $\Gamma$:
\begin{equation}
  \label{eq_full_l_system}
      \mathscr{L}(\Gamma') = \bigcup_{\Gamma ' \prec \widetilde\Gamma} \mathscr{L}_{\widetilde\Gamma}(\Gamma').
\end{equation}
\noindent
For instance, for $\Gamma = A_7$, we have:
\begin{equation*}
 \begin{split}
      & \mathscr{L}(A_7) = \mathscr{L}_{A_8}(A_7) \cup \mathscr{L}_{D_8}(A_7) \cup \mathscr{L}_{E_8}(A_7).
  \end{split}
\end{equation*}

If the diagram $\widetilde\Gamma$ in \eqref{eq_part_linkage_syst} is the $A$-type
(resp. $D$-type, resp. $E$-type) Dynkin diagram,
the component $\mathscr{L}_{\widetilde\Gamma}(\Gamma')$ is called the {\it $A$-component}
(resp. {\it $D$-component}, resp. {\it $E$-component}) in the full linkage system \eqref{eq_full_l_system}.
In {\bf Theorem \ref{th_size_Dtype}}, the $D$- and $E$-components,
the size of each component for
the full linkage systems $\mathscr{L}(D_l)$ and $\mathscr{L}(D_l(a_k))$ are found.
The results are shown in Table \ref{tab_number_linkages}.

\begin{table}[!ht]
\centering
\renewcommand{\arraystretch}{1.4}
  \begin{tabular} {|c|c|c|c|c|}
  \hline
     $\Gamma$ &  Number of  &
      \multicolumn{3}{c|}{Number of linkage diagrams,}  \cr
       &  components   & \multicolumn{3}{c|}{and $p = \mathscr{B}^{\vee}_{\Gamma}(\gamma^{\nabla})$}    \cr
      \cline{3-5}
         &    & $D$-components   & $E$-components &
          \hspace{1.8mm} In all \hspace{1.8mm} \cr
      \cline{3-4}
          &     & $p = 1$ & $p = \frac{l}{4}$ & \\
    \hline
       $D_4$, ~$D_4(a_1)$ &   $3$   & $3 \times 8 = 24$ & - &  $24$ \\
    \hline
       $D_5$, ~$D_5(a_1)$ &   $3$   &  $10$ & $2\times16 = 32$   & $42$ \\
    \hline
       $D_6$, $D_6(a_1)$, $D_6(a_2)$  &   $3$ & $12$  & $2\times32 = 64$  &    $76$ \\
    \hline
       $D_7$, $D_7(a_1)$, $D_7(a_2)$ &   $3$   & $14$ & $2\times64 = 128$  &  $142$ \\
    \hline
       $D_l$, ~$D_l(a_k)$, $l > 7$ &   $1$  & $2l$ & -  &  $2l$ \\
   \hline    
\end{tabular}
  ~\\ \vspace{2mm}
  \caption{Number of linkage diagrams in extensions for diagrams $D_l$ and $D_l(a_k)$}
  \label{tab_number_linkages}
\end{table}

For any pair of homogeneous Carter diagrams $\Gamma$ and  $\Gamma'$
the sizes of linkage systems $\mathscr{L}(\Gamma)$ and $\mathscr{L}(\Gamma')$
are the same ({\bf Proposition \ref{prop_sizes}}):
$\mid \mathscr{L}(\Gamma) \mid \;\; = \;\; \mid \mathscr{L}(\Gamma') \mid.$

The linkage systems $\mathscr{L}(D_4)$
(resp. $\mathscr{L}(D_4(a_1))$,  resp. $\mathscr{L}(D_l$) for $l >4$ , resp. $\mathscr{L}(D_l(a_k))$ for $l >4$)
are depicted in Fig.~\ref{D4_loctets} (resp. Fig.~\ref{D4a1_loctets}, resp. Fig.~\ref{Dlpu_linkages_2},
resp. Fig.~\ref{Dk_al_linkages}).

One of two orbits of the $E$-component of the linkage system $\mathscr{L}(D_5)$ are depicted in Fig.~\ref{D5pu_loctets}.
Two $E$-components of the linkage system $\mathscr{L}(D_6)$ are depicted in Fig.~\ref{D6pu_loctets}.

\subsection{The Weyl group of quadratic form $\mathscr{B}^{\vee}_{\Gamma}$}
  \label{sec_Weyl_gr}

The Weyl group of $\mathscr{B}^{\vee}_{\Gamma}$ is generated by dual reflections
$\{s^{*}_{\tau_i} \mid \tau_i \in S \}$, where
\begin{equation*}
       s^{*}_{\tau_i}{\gamma}^{\nabla} :=
             {\gamma}^{\nabla} - \langle {\gamma}^{\nabla}, \tau_i \rangle \tau_i^{\nabla}.
\end{equation*}

If ${\gamma}^{\nabla}$ is a linkage label vector then
values $(s^*_{\tau_i}{\gamma}^{\nabla})_{\tau_k}$ in \eqref{eq_dual_refl} belong to the set $\{-1, 0, 1\}$,
i.e., we get new linkage label vector $s^*_{\tau_i}{\gamma}^{\nabla}$ ({\bf Proposition \ref{prop_new_linkage}}).
The action of reflections $s_{\tau_i}$ and dual reflections $s^*_{\tau_i}$ related as follows:
\begin{equation*}
       (s_{\tau_i}\gamma)^{\nabla} = (s^*_{\tau_i}\gamma^{\nabla}).
\end{equation*}
The dual reflections $s^*_{\tau_i}$ preserve the inverse quadratic form $\mathscr{B}^{\vee}_{\Gamma}$
({\bf Proposition \ref{prop_link_diagr_conn}}):

\begin{equation*}
  \mathscr{B}^{\vee}_{\Gamma}(s^*_{\tau_i} {\gamma}^{\nabla}) =
     \mathscr{B}^{\vee}_{\Gamma}({\gamma}^{\nabla}) \text{ for } {\gamma}^{\nabla} \in L^{\nabla}.
\end{equation*}

The rational number $p = \mathscr{B}^{\vee}_{\Gamma}({\gamma}^{\nabla})$
is the invariant characterizing the given orbit, since
all elements of this orbit have the same value $\mathscr{B}^{\vee}_{\Gamma}({\gamma}^{\nabla})$,
see Table \ref{tab_number_linkages}.

\subsection{Loctets}

A spindle-like subset consisting of $8$ linkage label vectors is called a {\it loctet} (linkage octet).

\begin{figure}[!ht]
\centering
\includegraphics[scale=0.34]{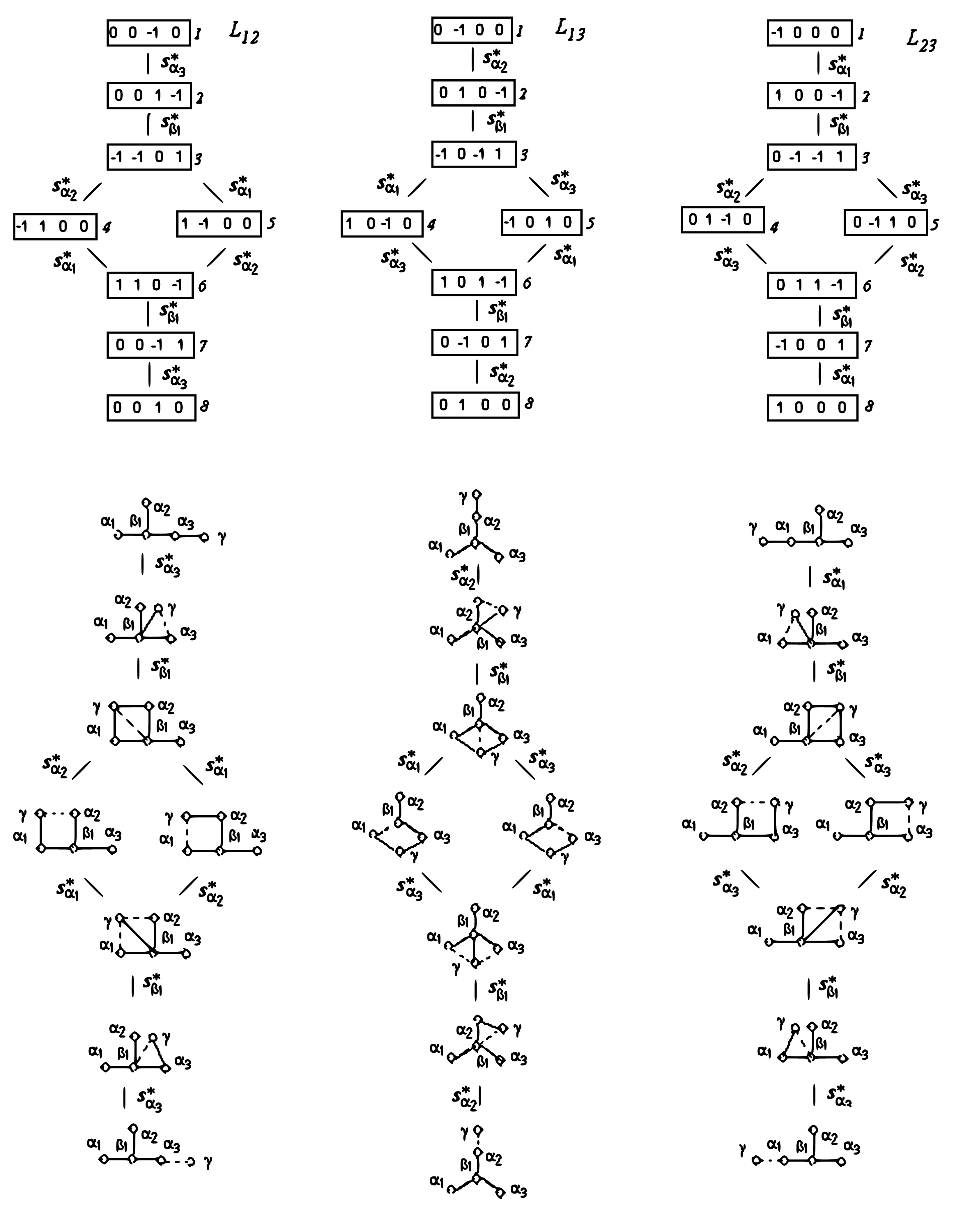}
 \caption{Three components of the linkage system $\mathscr{L}(D_4)$ contain $24$ elements.
 The linkage label vectors are shown above, the corresponding linkage diagrams are below.}
\label{D4_loctets}
\end{figure}

The $D$-component of $\mathscr{L}(D_4)$ consists of three orbits, each orbit consists of one loctet, Fig.~\ref{D4_loctets}.

\begin{figure}[!ht]
\centering
\includegraphics[scale=0.29]{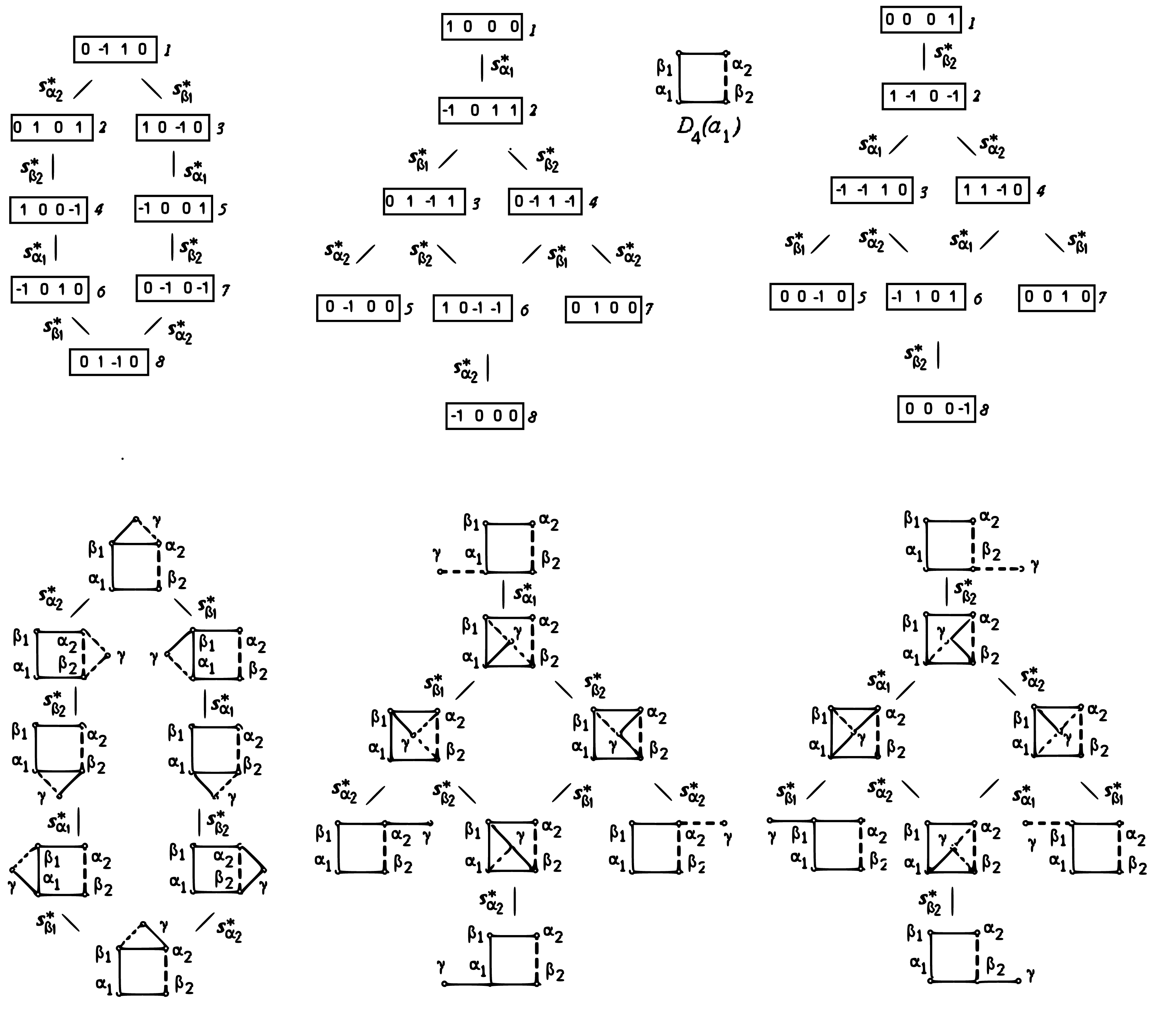}
 \caption{Three components of the linkage system $\mathscr{L}(D_4(a_1))$ contain $24$ elements.
 The linkage label vectors are shown above, the corresponding linkage diagrams are below.}
\label{D4a1_loctets}
\end{figure}

The $D$-component of $\mathscr{L}(D_5)$ consists of three orbits, two of which consist of one loctet, Fig.~\ref{D4a1_loctets}.

The $E$-component of $\mathscr{L}(D_5)$ consists of two orbits, each orbit consists of two loctets, Fig.~\ref{D5pu_loctets}.

The $E$-component of $\mathscr{L}(D_6)$ consists of two orbits, each orbit consists of four loctets, Fig.~\ref{D6pu_loctets}.

The $E$-component of $\mathscr{L}(D_7)$ consists of two orbits, each orbit consists of eight loctets, the corresponding
figures can be found in \cite[Figs. C.63, C.64]{St14}.

%% file: 3QuadrForm.tex
\section{\bf Quadratic forms associated with Carter diagrams}
  \label{sect_Cartan}

Let $L$ be the subspace spanned by the vectors $S = \{ \tau_1,\dots,\tau_n \}$.
We write this fact as follows:
\begin{equation*}
  L = [\tau_1,\dots,\tau_n],
\end{equation*}
The subspace $L$ is said to be the $S$-{\it associated subspace}.
Let $B$ be the Cartan matrix corresponding to the root system $\varPhi$.

 \begin{proposition}
   \label{prop_restr_forms}
 {\rm(i)}  The restriction of the bilinear form associated with the Cartan matrix
  $B$ on the subspace $L$ coincides with the bilinear form associated with the
  partial  Cartan matrix $B_{\Gamma}$, i.e., for any pair of vectors $v, u \in L$, we have
 \begin{equation}
   \label{restr_q}
       (v, u)_{\botG} = (v, u), \text{ and }
       \mathscr{B}_{\Gamma}(v) = \mathscr{B}(v).
 \end{equation}

 {\rm(ii)} For every Carter diagram, the matrix $B_{\botG}$ is positive definite.
 \end{proposition}
\PerfProof
 (i) From \eqref{canon_dec_2} we deduce:
 \begin{equation*}
     (v, u)_{\botG} =   (\sum\limits_i{t_i{\tau_i}}, \sum\limits_j{q_j{\tau_j}})_{\botG} =
      \sum\limits_{i,j}t_i{q}_j(\tau_i, \tau_j)_{\botG} = \sum\limits_{i,j}t_i{q}_j(\tau_i, \tau_j) =
      (v, u).
 \end{equation*}
 (ii) This follows from (i).
\qed
~\\

 If $\Gamma$ is a Dynkin diagram, the partial Cartan matrix $B_{\Gamma}$ is the Cartan matrix
 associated with $\Gamma$.  By \eqref{restr_q} the matrix $B_{\Gamma}$ is positive definite.
 The symmetric bilinear form associated with $B_{\Gamma}$ is denoted by $(\cdot, \cdot)_{\Gamma}$ and the
 corresponding quadratic form is denoted by $\mathscr{B}_{\Gamma}$.

 \subsection{The inverse quadratic form $\mathscr{B}^{\vee}_{\Gamma}$}
   \label{sec_inverse_qvadr}

\subsubsection{Action of $B^{-1}_{\Gamma}$ on the space $L$}
Let $S = \{\tau_1,\dots,\tau_l\}$ be a $\Gamma$-set, $L$ - the space spanned by roots of $S$,
 $\gamma$ be a root in $L$:  
\begin{equation*}
   \gamma = t_1\tau_1 + \dots + t_l\tau_l,
\end{equation*}
where $t_i$ are some rational numbers. Then, for $\gamma \in L$, we have
 \begin{equation}
  \label{eq_lin_depend2}
  \gamma^{\nabla} =
  \left (
    \begin{array}{c}
      (\gamma, \tau_1) \\
      \dots \\
      (\gamma, \tau_l) \\
    \end{array}
  \right ) =
  \left (
    \begin{array}{c}
      \sum t_i(\tau_i, \tau_1) \\
      \dots \\
      \sum t_i(\tau_i, \tau_l) \\
    \end{array}
  \right ) =
  B_{\Gamma}
  \left (
    \begin{array}{c}
      t_1 \\
      \dots\\
      t_l \\
    \end{array}
   \right ) = B_{\Gamma}\gamma,  \text{ and }
    \gamma = {B}^{-1}_{\Gamma}\gamma^{\nabla}.
 \end{equation}

\subsubsection{The projection of the linkage root}
  \label{sec_normal_mu}
Let $L^{\perp}$ be the orthogonal complement of $L$ to $L(\gamma)$
in the sense of the symmetric bilinear form $(\cdot, \cdot)$
associated with the root system $\varPhi$:
\begin{equation*}
   L(\gamma) =  L \oplus L^{\perp}.
\end{equation*}
Let $\gamma_L$ be the projection of the linkage root $\gamma$ on $L$. For
any root $\theta \in L(\gamma)$ such that $\theta \not\in L$, we
have $L(\theta) = L(\gamma)$, and $\theta$ is uniquely decomposed
into the following sum:
\begin{equation}
 \label{eq_theta_decomp}
   \theta =  \theta_L + \mu, \quad \text{ where } \quad \theta_L \in L, \quad \mu \in L^{\perp}.
\end{equation}
~\\
Given any vector $\theta$ by decomposition \eqref{eq_theta_decomp}, we introduce also the
{\it conjugate vector} $\overline{\theta}$ as follows:
\begin{equation*}
   \theta =  \theta_L + \mu, \quad  \overline{\theta} = \theta_L - \mu, \quad \text{ where } \quad \theta_L \in L, \quad \mu \in L^{\perp}.
\end{equation*}

 \begin{figure}[!ht]
\centering
\includegraphics[scale=0.2]{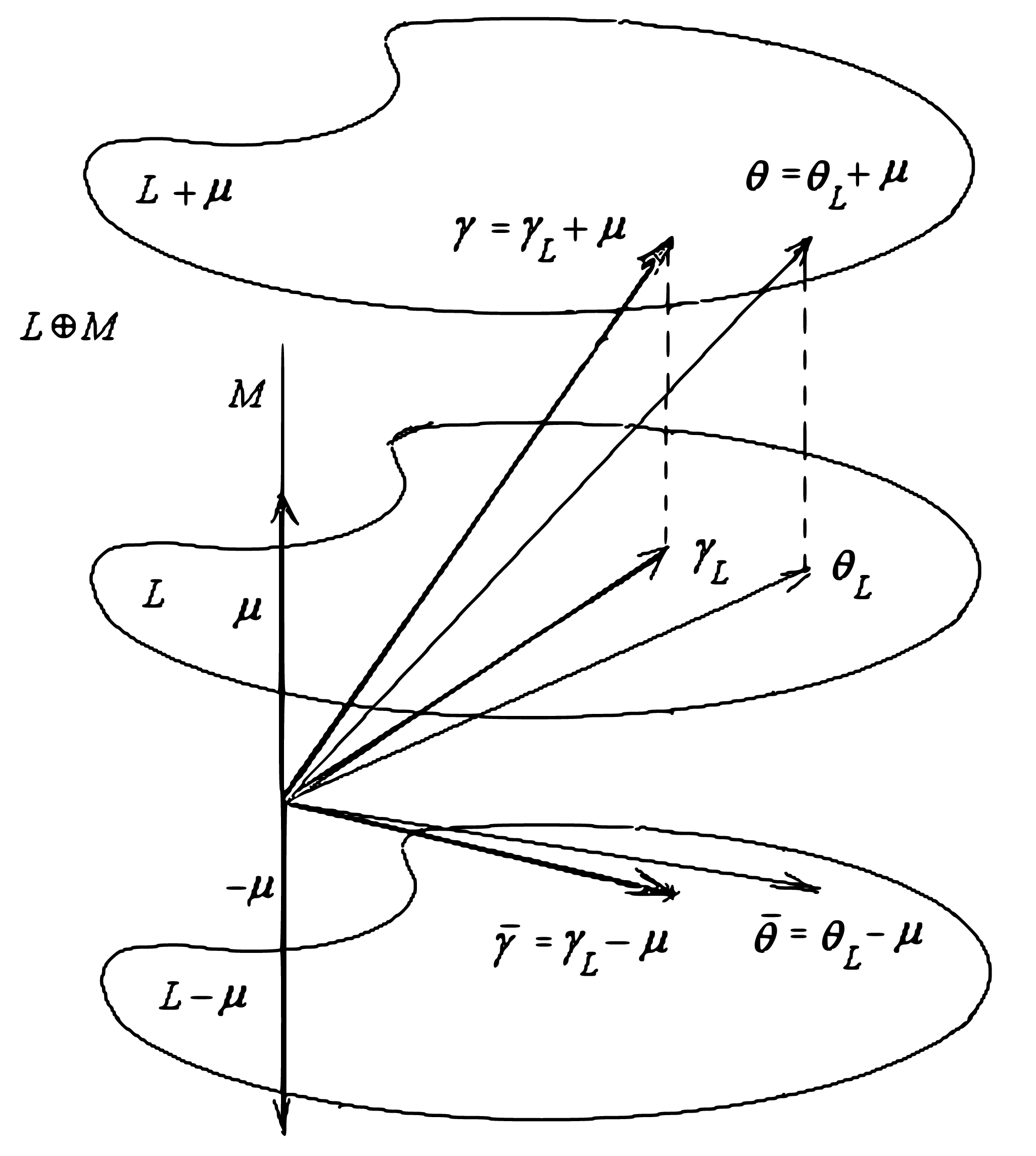}
 \caption{\hspace{3mm}The roots $\gamma = \gamma_L + \mu$ and $\overline{\gamma} = \gamma_L - \mu$.}
\label{L_and_mu}
\end{figure}
~\\
 The vector $\mu$ is said to be the {\it normal extending vector}, see Fig.~\ref{L_and_mu}.
 The following proposition collects a number of properties of
 the projection $\gamma_L$ of the linkage root $\gamma$, the linkage label vector $\gamma^{\nabla}$, and
 the normal extending vector $\mu$:

\begin{proposition}
  \label{prop_unique_val}
  \rm{(i)} The linkage label vector $\gamma^{\nabla}$ and the projection $\gamma_L$ are related as follows:
 \begin{equation}
  \label{eq_labels_proj}
   \gamma^{\nabla} = {\overline\gamma}^{\nabla} =  \gamma_L^{\nabla} = B_{\Gamma}\gamma_L.
 \end{equation}

  \rm{(ii)} The vector $\mu$ is, up to sign, a fixed vector for any root $\theta \in L(\gamma)\backslash L$.

  \rm{(iii)} The value $\mathscr{B}_{\Gamma}(\theta_L)$ is constant for any root $\theta \in L(\gamma)\backslash L$.

  \rm{(iv)} The vector $\theta = \theta_L + \mu$ is a root  if and only if
  the conjugate vector $\overline{\theta} = \theta_L - \mu$ is a root.
  In addition, $\theta \in L(\gamma)\backslash L$ if and only if
   $\overline{\theta} \in L(\overline\gamma)\backslash L$.

  \rm{(v)} Let $\theta$ is a root,  $\theta \not\in L$. Then, vector $\theta_L + t\mu$ is a root if and only if
  $t = \pm{1}$.

  \rm{(vi)} If $\delta$ is a root in $L(\gamma)\backslash L$ such that $\delta^{\nabla} = \gamma^{\nabla}$,
   then $\delta = \gamma$ or $\delta = \overline\gamma$.

  \rm{(vii)} For any root $\theta \in L(\gamma)\backslash L$, we have
\begin{equation*}
     \mathscr{B}_{\Gamma}^{\vee}(\theta^{\nabla}) = \mathscr{B}_{\Gamma}(\theta_L),
\end{equation*}
~\\
     and, $\mathscr{B}^{\vee}_{\Gamma}(\theta^{\nabla})$ is a constant for all roots  $\theta \in L(\gamma)\backslash L$.

 \end{proposition}

\PerfProof
 (i) Let $\gamma =  \gamma_L + \mu.$ Since $(\mu, \tau_i) = 0$ for any $i$, we have as follows:
 \begin{equation*}
      (\gamma, \tau_i) =  (\gamma_L, \tau_i) \; \text{ for any } i, \,\text{ i.e., }
      \gamma^{\nabla} = \gamma_L^{\nabla}.
 \end{equation*}
  Then, by \eqref{eq_lin_depend2} we have
 \begin{equation*}
   \gamma^{\nabla} = {\overline\gamma}^{\nabla} =  \gamma_L^{\nabla} = B_{\Gamma}\gamma_L.
 \end{equation*}

 (ii) Any root $\theta \in L(\gamma)$ such that $\theta \not\in L$ decomposes as $\tau + \gamma$
 for some $\tau \in L$. Since $\gamma = \gamma_L + \mu$, where $\gamma_L \in L$,  we have
 $\theta = \tau + \gamma_L + \mu$. Here, $\tau + \gamma_L \in L$. Let us put $\theta_L = \tau + \gamma_L$.
 Then $\theta = \theta_L + \mu$.
 ~\\

 (iii) Let $\mathscr{B}$ be the quadratic form associated with the root system $\varPhi$.
  By \eqref{eq_theta_decomp} we have $\theta_L \perp \mu$, and
  $\mathscr{B}(\theta) = \mathscr{B}(\theta_L) + \mathscr{B}(\mu)$.
  Here, $\mathscr{B}(\theta) = 2$ since $\theta$ is the root, and by (ii),
  we have $\mathscr{B}(\theta_L) = 2 - \mathscr{B}(\mu)$, i.e, $\mathscr{B}(\theta_L)$ is constant.
  By \eqref{restr_q}, we have $\mathscr{B}_{\Gamma}(\theta_L) = \mathscr{B}(\theta_L)$,
  i.e.,  $\mathscr{B}_{\Gamma}(\theta_L)$ is also constant for all $\theta \in L(\gamma)$.
~\\

 (iv) Let $\theta$ be a root, i.e., $\mathscr{B}(\theta) = \mathscr{B}(\theta_L) + \mathscr{B}(\mu) = 2$.
  Then for $\overline{\theta}$, we have
  $\mathscr{B}(\overline{\theta}) = \mathscr{B}(\theta_L) + \mathscr{B}(-\mu) = 2$ as well.
  By \eqref{eq_Kac} $\overline{\theta}$ is a root as well.
  Further, $\theta \in L(\gamma)$ means that $\theta = \gamma + \tau$, where $\tau \in L$.
  Thus,
   $\theta = \gamma_L + \mu + \tau$. Put  $\theta_L =  \gamma_L + \tau$. Then $\theta_L \in L$,  and
\begin{equation*}
  \theta = \theta_L + \mu, \;\; \overline\theta = \theta_L - \mu  = \gamma_L + \tau - \mu =
  \overline\gamma + \tau \in L(\overline\gamma).
\end{equation*}

 (v) Consider root $\theta = \theta_L \pm\mu$. Since $\theta_L \perp \mu$,  the following relations hold:
\begin{equation*}
   \mathscr{B}(\theta_L + t\mu) = \mathscr{B}(\theta \pm\mu + t\mu) = \mathscr{B}(\theta +  (t \pm{1})\mu) =
   \mathscr{B}(\theta) + (t \pm{1})^2\mathscr{B}(\mu).
\end{equation*}
   Because $\theta$ is a root, then
   $\mathscr{B}(\theta_L + t\mu) = \mathscr{B}(\theta) = 2$ and $t = \pm{1}$.
~\\

 (vi) By  $\delta^{\nabla} = \gamma^{\nabla}$ and \eqref{eq_labels_proj}, the roots $\delta$ and $\gamma$
 have the same projections on $L$: $\delta_L = \gamma_L$.
 If $\delta = \delta_L + t\mu = \gamma_L + t\mu$ is a root, then by (v) we have $t =\pm{1}$.
 Thus, $\delta = \gamma$ or $\overline{\gamma}$.
~\\

 (vii) By heading (i), since $\theta_L \in L$, we have $\theta_L^{\nabla} = B_{\Gamma}\theta_L$.  Thus,
 \begin{equation*}
      \mathscr{B}^{\vee}_{\Gamma}(\theta^{\nabla}) =
      \langle B_{\Gamma}^{-1}\theta^{\nabla}, \theta^{\nabla} \rangle =
     \langle {\theta_L},  B_{\Gamma}\theta_L \rangle = \mathscr{B}_{\Gamma}(\theta_L). \\
 \end{equation*}
\qed

\subsection{The linkage root theorem}
   \index{linkage root}
   \index{linkage root criterion}

 In this section we give a criterion for a given vector to be a linkage root.

\begin{theorem}[Linkage root criterion]
   \label{th_B_less_2}
     {\rm(i)} Let $\theta^{\nabla}$ be the linkage label vector
     corresponding to a certain root $\theta \in L(\gamma)$, i.e., $\theta^{\nabla} = B_{\Gamma}\theta_L$.
     The root $\theta$ is a linkage root,
     (i.e., $\theta$ is linearly independent of roots of $L$) if and only if
     \begin{equation}
       \label{eq_p_less_2}
          \mathscr{B}^{\vee}_{\Gamma}(\theta^{\nabla}) < 2.
     \end{equation}

    {\rm(ii)} Let $\theta \in L(\gamma)$ be a root connected only
     with one $\tau_i \in S$. The root $\theta$ is a linkage root if and only if
     \begin{equation}
       \label{eq_bii_less_2}
          b^{\vee}_{ii} < 2,
     \end{equation}
     where $b^{\vee}_{ii}$ is the $i$th diagonal element of $B^{-1}_{\Gamma}$.
\end{theorem}

\PerfProof
    (i) By  Proposition \ref{prop_unique_val}(vi) and  Proposition \ref{prop_restr_forms} we have
\begin{equation}
    \mathscr{B}^{\vee}_{\Gamma}(\theta^{\nabla}) = \mathscr{B}_{\Gamma}(\theta_L) = \mathscr{B}(\theta_L) =
    \mathscr{B}(\theta) - \mathscr{B}(\mu) \leq 2.
\end{equation}
  Then, $\mathscr{B}^{\vee}_{\Gamma}(\theta^{\nabla}) = 2$ if and only if $\mathscr{B}(\mu) = 0$.
  If $\mathscr{B}^{\vee}_{\Gamma}(\theta^{\nabla}) = 2$ then $\mathscr{B}(\mu) = 0$,  $\mu = 0$
  and $\theta = \theta_L$. Therefore, $\theta$ is linearly depends on vectors of $L$.
  If $\mathscr{B}^{\vee}_{\Gamma}(\theta^{\nabla}) < 2$ then $\mathscr{B}(\mu) > 0$, i.e., $\mu \neq 0$
  and $\theta$ is linearly independent of roots of $L$.

    (ii) We have
     \begin{equation*}
        \theta^{\nabla} \quad = \quad
                \left (
                 \begin{array}{c}
                    (\theta, \tau_1) \\
                    \dots \\
                    (\theta, \tau_i) \\
                    \dots \\
                    (\theta, \tau_l) \\
                 \end{array}
                  \right )
                \quad = \quad
                 \left (
                 \begin{array}{c}
                    0 \\
                    \dots \\
                    \pm{1} \\
                    \dots \\
                    0 \\
                 \end{array}
                  \right ),
     \end{equation*}
  and $\mathscr{B}^{\vee}_{\Gamma}(\theta^{\nabla}) = b^{\vee}_{ii}$.
  Thus,  statement (ii) follows from (i).
  \qed
\bigskip

\subsection{Quadratic forms for homogeneous Carter diagrams}

\begin{proposition}
  \label{prop_quadr_1}
   Let $\{\Gamma', \Gamma\}$ be a pair of Carter diagrams from the adjacency list \eqref{eq_pairs_trans},
   $M$ be the transition matrix transforming $\Gamma'$-set $S'$ to $\Gamma$-set $S$ as in \S\ref{sec_trans_mappings}.
   The vector $\gamma^{\nabla}$ is the linkage label vector for $\Gamma$ if and only of
   \, ${}^t{M}\gamma^{\nabla}$ is the linkage label vector for $\Gamma'$.
\end{proposition}

\PerfProof
By Theorem \ref{th_invol}, the transition matrix $M$ maps the $\Gamma'$-set $S'$ to the $\Gamma$-set $S$:
\begin{equation*}
   M\tau'_i = \tau_i, \; \text{ where } \; \tau'_i \in S', \tau_i \in S.
\end{equation*}
Let $B_{\Gamma'}$ and $B_{\Gamma}$ be the partial Cartan matrices for $\Gamma'$ and $\Gamma$.
The transition matrix $M$ transforms $S'$ to $S$ and relates the partial Cartan matrices
$B_{\Gamma'}$ and $B_{\Gamma}$ as follows:
\begin{equation}
 \label{eq_relate_bases}
   {}^t{M} \cdot B_{\Gamma'} \cdot M = B_{\Gamma}, \\
\end{equation}
see \cite[\S4.7]{St23}.
Since $M$ is an involution, we get
\begin{equation*}
  M \cdot B^{-1}_{\Gamma'} \cdot {}^tM = B^{-1}_{\Gamma}.
\end{equation*}
Further,
\begin{equation*}
 \begin{split}
   \mathscr{B}^{\vee}_{\Gamma}(\gamma^{\nabla}) = & \langle B^{-1}_{\Gamma}\gamma^{\nabla}, \gamma^{\nabla} \rangle =
     \langle  B^{-1}_{\Gamma'} \cdot {}^tM\gamma^\nabla, {}^tM\gamma^{\nabla} \rangle =
      \mathscr{B}^{\vee}_{\Gamma'}({}^tM\gamma^{\nabla}).
 \end{split}
\end{equation*}
Then, $ \mathscr{B}^{\vee}_{\Gamma}(\gamma^{\nabla}) < 2$ \; if and only if \;
$\mathscr{B}^{\vee}_{\Gamma'}({}^tM\gamma^{\nabla}) < 2$. By Theorem \ref{th_B_less_2},
we get what we need. \qed
~\\

\subsection{The chain of homogeneous pairs}
 \label{sec_chain}
Consider any homogeneous pair of Carter diagrams $\{\Gamma', \Gamma\}$.
There exists the chain of homogeneous pairs from the adjacency list \eqref{eq_pairs_trans}
which connect $\Gamma'$ and $\Gamma$ as follows:
\begin{equation}
  \label{eq_chain_1}
   \{\{\Gamma_0,\Gamma_1\},\{\Gamma_1,\Gamma_2\},\dots, \{\Gamma_{k-1},\Gamma_k \}, \text{ where }
     \Gamma_0 = \Gamma' \text{ and } \Gamma_k =\Gamma. \\
\end{equation}
Let $S_j = \{\tau^{j}_1, \dots, \tau^{j}_n\}$  be the associated $\Gamma_j$-set,
$M^{(j)}$ be the transition matrices connecting $j$th pair in \eqref{eq_chain_1}, where $j = 0,\dots k-1$,
see Theorem \ref{th_invol}.
Then,
\begin{equation*}
    \begin{cases}
      M^{(j)}(\tau^j_i) = \tau^{j+1}_i  \text{ for some } i \in \{1,\dots,n\}  , \\
      M^{(j)}(\tau^j_k) = \tau^{j}_k  \quad \text{ for } k \neq i. \\
   \end{cases}   \\
\end{equation*}
Two Carter diagrams obtained from each other by a sequence of reflections \eqref{eq_similarity_0},
are called {\it similar Carter diagrams}, see \cite[Fig.~4]{St23}.
In fact, mappings $L_{\tau_i}$ can be considered as a degenerated case
of transition matrices $M$ of Theorem \ref{th_invol}.
The Carter diagrams in the adjacency list \eqref{eq_pairs_trans} are given up
to similarity, so using the transition matrices $M^{(j)}$ together with mappings $L_{\tau_i}$
we can construct the map $F$ transforming the partial Cartan matrices $B_{\Gamma'}$
and $B_{\Gamma}$ to each other, also the linkage roots for $\Gamma'$ and $\Gamma$
to each other, see \cite[\S4.2.1]{St23}:
\begin{equation}
\label{eq_matr_F}
 \begin{split}
  & F\tau'_i = \tau_i, \quad
   {}^t{F} \cdot B_{\Gamma'} \cdot F = B_{\Gamma}, \quad
    \mathscr{B}^{\vee}_{\Gamma}(\gamma^{\nabla}) = \mathscr{B}^{\vee}_{\Gamma'}({}^tF\gamma^{\nabla}).
 \end{split}
\end{equation}

\subsection{Spectrum of the Dynkin and Carter diagram}
The spectrum of a Dynkin diagram (resp. Carter diagram) is defined as
the spectrum of the corresponding Cartan matrix (resp. partial Cartan matrix).

\begin{proposition}
  \label{prop_spectr}
  The spectrum of any Dynkin or Carter diagram lies in the open interval $(0, 4)$.
\end{proposition}
\PerfProof
Note that for any transition matrix $M$ and any similarity matrix $L_{\tau}$
the determinant is equal to $-1$, see \eqref{th_transition}.
Then, by \eqref{eq_relate_bases}  and \eqref{eq_matr_F}, we have
\begin{equation}
  \label{eq_determ}
   \det B_{\Gamma'} = \det B_{\Gamma}
\end{equation}
for any homogeneous pair of Carter diagrams $\{\Gamma', \Gamma\}$.
Let $\Gamma$ be the Dynkin diagram from the corresponding homogeneous class.
By \eqref{eq_determ} the spectrum of any Carter diagram
coincides with the spectrum of $\Gamma$.

Now consider the spectrum of Dynkin diagrams.
By \cite[eq.(3.2)]{St08} and \cite[Prop.3.2]{St08} eigenvalues $\lambda$ of the Coxeter transformation
and the eigenvalues $\rho$ of $B_{\Gamma}$ are related as follows:
\begin{equation}
  \label{eq_Coxeter}
  \begin{split}
  &\frac{(\lambda + 1)^2}{4\lambda} = (\frac{\rho}{2} - 1)^2 \quad \text{ or } \\
  &      \lambda + 2 + \frac{1}{\lambda} = (\rho - 2)^2.
  \end{split}
\end{equation}

Since eigenvalues of the Coxeter transformation are roots of unity, see \cite{Bo02},
any $\lambda$ is equal to  $\cos\varphi + i\sin\varphi$ for some angle $\varphi$.
By \eqref{eq_Coxeter}
\begin{equation*}
  2\cos\varphi + 2 \;=\; (\rho - 2)^2, \; \text{ i.e., } \; (\rho-2)^2 \;=\; 4\cos^2\frac{\varphi}{2} \;\leq\; 4.
\end{equation*}

Therefore, $0 \le \rho \le 4$.
By Proposition \ref{prop_restr_forms}(ii) $B_{\Gamma}$ is positive definite, then $\rho > 0$.
For the finite Weyl groups, $1$ is not an eigenvalue of the Coxeter transformation, then $\rho \neq 4$.
Thus, $0   <  \rho  <  4$.
\qed
~\\

Using the spectrum of the Cartan matrix could be an invariant description of the Dynkin and Carter diagrams,
see \cite[\S1.5.3]{St23}, \cite{MS07}.
~\\

%% file: 4WeylGr_Dual.tex
\subsection{The Weyl group of the quadratic form}

\subsubsection{Dual reflections}
    \label{sec_dual Weyl_gr}
 Let $\Gamma$ be a Carter diagram,
$S = \{\tau_1,\dots,\tau_l\}$ be some
$\Gamma$-set, and $L$ be $S$-associated subspace:
\begin{equation*}
  L = [\tau_1,\dots,\tau_l].
\end{equation*}
 For any  $\tau_i \in S$, we define vectors $\tau_i^{\nabla}$ as follows:
   \begin{equation}
      \label{eq_tau_vee}
        \tau_i^{\nabla} := B_{\Gamma}{\tau_i}.
   \end{equation}
 Let ${L}^{\nabla}$ be the linear space spanned by the vectors $\tau_i^{\nabla}$, where $\tau_i \in S$.
   The map $\tau_i \longrightarrow \tau_i^{\nabla}$ given by \eqref{eq_tau_vee}
   is expanded to the linear mapping $L \longrightarrow {L}^{\nabla}$, and
   \begin{equation}
      \label{eq_gamma_vee}
       u^{\nabla} = B_{\Gamma}{u} =
       \left (
       \begin{array}{c}
         (u, \tau_1) \\
          \dots \\
         (u, \tau_l) \\
       \end{array}
       \right ) \text{ for any } u \in L,
   \end{equation}
   see \eqref{eq_lin_depend2}.
   Consider the restriction of the reflection $s_{\tau_i}$ on the subspace $L$.
   For any $v \in L$, by Proposition \ref{prop_restr_forms} we have:
 \index{dual reflection $s^{*}_{\tau_i}$}
 \index{$s^{*}_{\tau_i}$, dual reflection}
   \begin{equation}
     \label{refl_tau}
       s_{\tau_i}v = v  - 2\frac{(\tau_i, v)}{(\tau_i, \tau_i)}\tau_i =
                      v  - (\tau_i, v)_{\botG}\tau_i =
                      v  - \langle B_{\Gamma}\tau_i, v \rangle \tau_i =
                      v - \langle \tau_i^\nabla, v \rangle \tau_i.
   \end{equation}
   We define the {\it dual reflection} $s^{*}_{\tau_i}$ acting on a vector $u^{\nabla} \in L^{\nabla}$ as follows:
   \begin{equation}
     \label{refl_tau_dual}
       s^{*}_{\tau_i}u^{\nabla} := u^{\nabla} - \langle u^{\nabla}, \tau_i \rangle \tau_i^{\nabla}.
   \end{equation}
   Let $W_S$ (resp. $W_S^{\vee}$) be the group generated by reflections
   $\{s_{\tau_i} \mid \tau_i \in S \}$ (resp. $\{s^{*}_{\tau_i} \mid \tau_i \in S \})$,
   where $S = \{\tau_1,\dots,\tau_l\}$.  The group $W_S$ is uniquely determined by the quadratic
   form $\mathscr{B}_{\Gamma}$.
   In Proposition \ref{prop_link_diagr_conn} it will be shown that $W_S$ (resp. $W_S^{\vee}$) preserves
   $\mathscr{B}_{\Gamma}$ (resp. $\mathscr{B}^{\vee}_{\Gamma}$), so $W_S$ (resp. $W_S^{\vee}$)
   is called the {\it Weyl group of} $\mathscr{B}_{\Gamma}$ (resp. $\mathscr{B}^{\vee}_{\Gamma}$)\footnotemark[1].
   \footnotetext[1]{
   The Weyl group of the integral quadratic form and their roots
   were studied by A.~Roiter in \cite{R78}, and then this study was continued in
   \cite{GR97}, \cite{S18}, \cite{MM19}, \cite{MZ22}  and others.}

   \begin{proposition}
    \label{prop_contagr_1}
     {\rm(i)} For any $\tau_i \in S$, we have
   \begin{equation}
      \label{refl_transp_2}
        s^{*}_{\tau_i} = {}^t{s}_{\tau_i} = {}^t{s}_{\tau_i}^{-1}.
   \end{equation}
     {\rm(ii)} The mapping
   \begin{equation*}
        \pi:  w \rightarrow {}^t{w}^{-1}
   \end{equation*}
     determines an isomorphism between $W_S$ and $W_S^{\vee}$.
   \end{proposition}

   \PerfProof
  (i) By \eqref{refl_tau} and \eqref{refl_tau_dual}, for any $v \in L, u^{\nabla} \in L^{\nabla}$, we have:
   \begin{equation*}
     \begin{split}
      & \langle  s^{*}_{\tau_i}u^{\nabla}, v \rangle =
       \langle  u^{\nabla} - \langle u^{\nabla}, \tau_i \rangle \tau_i^{\nabla}, v \rangle =
        \langle u^{\nabla}, v \rangle - \langle u^{\nabla}, \tau_i \rangle \langle v, \tau_i^{\nabla} \rangle, \\
      &   \langle  u^{\nabla}, s_{\tau_i}v \rangle =
        \langle  u^{\nabla} , v - \langle \tau_i^{\nabla}, v \rangle \tau_i \rangle =
        \langle u^{\nabla}, v \rangle - \langle \tau_i^{\nabla}, v \rangle \langle u^{\nabla}, \tau_i \rangle.
     \end{split}
   \end{equation*}
   Thus,
   \begin{equation*}
      \langle  s^{*}_{\tau_i}u^{\nabla}, v \rangle = \langle  u^{\nabla}, s_{\tau_i}v \rangle,
      \quad \text{ for any } v \in L, u^{\nabla} \in L^{\nabla},
   \end{equation*}
   and \eqref{refl_transp_2} holds.
~\\

   (ii) Let $\pi(w_i) = {}^t w^{-1}_i$ for $i = 1,2$. Then,
   \begin{equation*}
      \pi(w_1)\,\pi(w_2) = {}^t w^{-1}_1 \; {}^t w^{-1}_2 = {}^t (w^{-1}_2 \; w^{-1}_1) =
      {}^t (w_1 \; w_2)^{-1} = \pi(w_1 \; w_2).
   \end{equation*}
  \qed
 ~\\
 Let us put
   \begin{equation*}
      w^* := {}^t{w}^{-1}.
   \end{equation*}
 Then, $\pi(w) =  w^*$, and
   \begin{equation}
     \label{eq_isom}
    w^*_1 \; w^*_2 = (w_1 w_2)^*.
   \end{equation}
\smallskip

\subsubsection{Action of the Weyl group $W_S^{\vee}$}
   \label{sec_partial_groups}

  By \eqref{eq_tau_vee}, \eqref{eq_gamma_vee} and \eqref{refl_tau_dual}, for any $u^{\nabla} \in L^{\nabla}$,
  we have
   \begin{equation}
     \label{refl_tau_dual_2}
    \begin{split}
       &  s^{*}_{\tau_i}u^{\nabla} = u^{\nabla} - \langle u^{\nabla}, \tau_i \rangle \tau^{\nabla}_i =
             u^{\nabla} - u^{\nabla}_{\tau_i} (\tau^{\nabla}_i) =
             u^{\nabla} - u^{\nabla}_{\tau_i} (B_{\Gamma}\tau_i),
          \\
       & (s^{*}_{\tau_i}u^{\nabla})_{\tau_k} = u^{\nabla}_{\tau_k} - u^{\nabla}_{\tau_i}(B_{\Gamma}\tau_i)_{\tau_k}
                                  = u^{\nabla}_{\tau_k} - u^{\nabla}_{\tau_i}(\tau_i, \tau_k).
    \end{split}
   \end{equation}

  Let $u$ be a linkage root, then vector $u^{\nabla}$ is the linkage label vector, the coordinates of
  $u^{\nabla}$ belong to the set $\{-1, 0, 1\}$.
  By \eqref{refl_tau_dual_2} we get
    \begin{equation}
       \label{eq_dual_refl}
       (s^*_{\tau_i}u^{\nabla})_{\tau_k} =
        \begin{cases}
            -u^{\nabla}_{\tau_i}, & \text{ if } k = i, \\
            u^{\nabla}_{\tau_k} + u^{\nabla}_{\tau_i}, & \text{ if } \{\tau_k, \tau_i \}
            \text{ is a {\it solid} edge, i.e., } (\tau_k, \tau_i) = -1, \\
            u^{\nabla}_{\tau_k} - u^{\nabla}_{\tau_i}, & \text{ if }  \{\tau_k, \tau_i \}
            \text{ is a {\it dotted} edge, i.e., } (\tau_k, \tau_i) = 1, \\
            u^{\nabla}_{\tau_i}, & \text{ if } \tau_k \text{ and } \tau_i
             \text{ are not connected, i.e., } (\tau_k, \tau_i) = 0.
        \end{cases}
    \end{equation}

\begin{proposition}[Action of dual reflection]
   \label{prop_new_linkage}
    If $u^{\nabla}$ is a linkage label vector then
    values $(s^*_{\tau_i}u^{\nabla})_{\tau_k}$ in \eqref{eq_dual_refl} belong to the set $\{-1, 0, 1\}$,
    i.e., we obtain a new linkage label vector $s^*_{\tau_i}u^{\nabla}$.
 \end{proposition}
 \PerfProof
    Let $\{\tau_k, \tau_i\}$ be a solid edge. If $u^{\nabla}_{\tau_k} = u^{\nabla}_{\tau_i} = 1$
    (resp. $u^{\nabla}_{\tau_k} = u^{\nabla}_{\tau_i} = -1$) then roots $\{-u, \tau_k, \tau_i\}$
    (resp. $\{u, \tau_k, \tau_i\}$) constitute the root system corresponding to the
    extended Dynkin diagram $\widetilde{A}_2$, which is impossible,
    see \cite[Lemma A.1]{St17}.
    For remaining pairs $\{u^{\nabla}_{\tau_k}, u^{\nabla}_{\tau_i}\}$, we have
    $-1 \leq u^{\nabla}_{\tau_k} + u^{\nabla}_{\tau_i} \leq 1$. Now, let $\{\tau_k, \tau_i\}$ be a dotted edge.
    If $u^{\nabla}_{\tau_k} = 1$ and $u^{\nabla}_{\tau_i} = -1$ (resp. $u^{\nabla}_{\tau_k} = -1$ and $u^{\nabla}_{\tau_i} = 1$)
    then roots $\{u, -\tau_k, \tau_i\}$ (resp. $\{u, \tau_k, -\tau_i\}$)
    constitute the root system $\widetilde{A}_2$, which is impossible.
    For remaining pairs $\{u^{\nabla}_{\tau_k}, u^{\nabla}_{\tau_i}\}$, we have
    $-1 \leq u^{\nabla}_{\tau_k} - u^{\nabla}_{\tau_i} \leq 1$.
\qed

\begin{remark}
  \label{rem_action}
{\rm
  Eq. \eqref{eq_dual_refl} is very useful for calculating linkage systems.
  Note that action of dual reflections in \eqref{eq_dual_refl} $s^*_{\tau_k}$ on the linkage label vector
  $u^{\nabla}$ is non-trivial if and only if $u^{\nabla}_k \ne 0$.}
\end{remark}

\begin{proposition}
  \label{prop_link_diagr_conn}
  {\rm(i)} For dual reflections $s^{*}_{\tau_i}$, the following relations hold:
  \begin{equation}
      \label{conn_dual_1}
        B_{\Gamma} s_{\tau_i} = s^*_{\tau_i} B_{\Gamma},
     \end{equation}
  \begin{equation}
     \label{conn_L_2}
    (s_{\tau_i}\gamma)^{\nabla} = s^{*}_{\tau_i}B_{\Gamma}\gamma_L = s^{*}_{\tau_i}\gamma^{\nabla}.
  \end{equation}
  {\rm(ii)} For $w^* \in W^{\vee}_S$, we have 
    \begin{equation*}
        (w\gamma)^{\nabla} = w^*{\gamma}^{\nabla}.
     \end{equation*}
  {\rm(iii)} The following relations hold:
    \begin{equation}
      \label{conn_L_5}
        \begin{split}
        & \mathscr{B}_{\Gamma} (s_{\tau_i} v)  =  \mathscr{B}_{\Gamma}(v)  \text{ for any } v \in L, \\
        & \mathscr{B}^{\vee}_{\Gamma}(s^*_{\tau_i} u^{\nabla}) =  \mathscr{B}^{\vee}_{\Gamma}(u^{\nabla}) \text{ for any } u^{\nabla} \in L^{\nabla}.
        \end{split}
     \end{equation}
     \end{proposition}

\PerfProof (i)  The equality \eqref{conn_dual_1} is true since for any $u, v \in L$ the
following is true:
 \begin{equation*}
    \begin{split}
     & (s_{\tau_i}{u}, v)_{\Gamma} = ({u}, s_{\tau_i}{v})_{\Gamma}, \text{ i.e., }
     \langle B_{\Gamma} s_{\tau_i}{u}, v \rangle = \langle B_{\Gamma} {u}, s_{\tau_i}{v} \rangle =
     \langle s^*_{\tau_i}{B}_{\Gamma} {u}, v \rangle \text{, and } \\
     & \langle (B_{\Gamma} s_{\tau_i} - s^*_{\tau_i}{B}_{\Gamma}) {u}, v \rangle = 0.
    \end{split}
 \end{equation*}
~\\
Let us consider eq.~\eqref{conn_L_2}. Since $(\tau_i, \mu) = 0$ for
any $\tau_i \in S$, and $s_{\tau_i}\mu = \mu$, we have
 \begin{equation*}
  \begin{split}
    (s_{\tau_i}\gamma)^{\nabla} =
      & \left (
      \begin{array}{c}
         (s_{\tau_i}\gamma, \tau_1) \\
         \dots \\
         (s_{\tau_i}\gamma, \tau_l) \\
      \end{array}
      \right ) =
      \left (
      \begin{array}{c}
         (s_{\tau_i}\gamma_L + \mu, \tau_1) \\
         \dots \\
         (s_{\tau_i}\gamma_L + \mu, \tau_l) \\
      \end{array}
      \right ) =
      \left (
      \begin{array}{c}
         (s_{\tau_i}\gamma_L, \tau_1) \\
         \dots \\
         (s_{\tau_i}\gamma_L, \tau_l) \\
      \end{array}
      \right ) =
      (s_{\tau_i}\gamma_L)^{\nabla}.
      \\
     \end{split}
  \end{equation*}
  Then, by \eqref{conn_dual_1} and \eqref{eq_labels_proj}
  \begin{equation*}
    (s_{\tau_i}\gamma)^{\nabla} =
    (s_{\tau_i}\gamma_L)^{\nabla} = B_{\Gamma}{s}_{\tau_i}\gamma_L =
    s^{*}_{\tau_i}B_{\Gamma}\gamma_L = s^{*}_{\tau_i}\gamma^{\nabla}.
  \end{equation*}

(ii)  Let $w = s_{\tau_1}s_{\tau_2}\dots{s}_{\tau_m}$ be a
decomposition of $w \in W$. Since  $s^*_{\tau} = {}^ts^{-1}_{\tau} =
{}^ts_{\tau}$,
 we deduce from \eqref{eq_isom} and \eqref{conn_L_2} the following:
\begin{equation*}
  \begin{split}
   & (w\gamma)^{\nabla} = (s_{\tau_1}s_{\tau_2}\dots{s}_{\tau_m}\gamma)^{\nabla} =
    s^*_{\tau_1}(s_{\tau_2}\dots{s}_{\tau_m}\gamma)^{\nabla} =
    s^*_{\tau_1}s^*_{\tau_2}({s}_{\tau_3}\dots{s}_{\tau_m}\gamma)^{\nabla} = \dots = \\
   & s^*_{\tau_1}{s}^*_{\tau_2}\dots{s}^*_{\tau_m}{\gamma}^{\nabla} =
     {w}^*{\gamma}^{\nabla}.
  \end{split}
\end{equation*}

 (iii)
 Further,  by \eqref{conn_dual_1} we have
 \begin{equation*}
    \begin{split}
    & \mathscr{B}_{\Gamma}(s_{\tau_i}{v}) = \langle B_{\Gamma} s_{\tau_i}{v}, s_{\tau_i}{v} \rangle =
     \langle s^*_{\tau_i}{B}_{\Gamma} {v}, s_{\tau_i}{v} \rangle =
     \langle {B}_{\Gamma} {v}, {v} \rangle = \mathscr{B}_{\Gamma}(v). \\
    \hspace{2.3cm} &
     \mathscr{B}_{\Gamma}^{\vee}(s^{*}_{\tau_i}{u^{\nabla}}) =
     \langle B_{\Gamma}^{\vee} s^{*}_{\tau_i}{u^{\nabla}}, s^{*}_{\tau_i}{u^{\nabla}} \rangle =
     \langle s_{\tau_i}{B}_{\Gamma}^{\vee} {u^{\nabla}}, s^{*}_{\tau_i}{u^{\nabla}} \rangle =
     \langle {B}_{\Gamma}^{\vee} {u^{\nabla}}, {u^{\nabla}} \rangle = \mathscr{B}_{\Gamma}^{\vee}(u^{\nabla}).
     \hspace{2.3cm}  \\
    \end{split}
 \end{equation*}
 so \eqref{conn_L_5} holds. \qed
 ~\\

  By Propositions \ref{prop_new_linkage} each linkage label vector belongs to some orbit
  of the group $W^{\vee}_S$ action.
  By Proposition \ref{prop_link_diagr_conn}(iii)
  all elements of some orbit have the same value $\mathscr{B}^{\vee}_{\Gamma}(u^{\nabla})$.
  The rational number $p = \mathscr{B}^{\vee}_{\Gamma}(u^{\nabla})$ is the invariant
  characterizing the given orbit.

  To find the preimage root $\gamma_1$ for
  some label vector $\gamma_1^{\nabla}$  is suffices to know the preimage for one element $\gamma$ on the orbit.
  By Proposition \ref{prop_link_diagr_conn}(ii) if $\gamma_1^{\nabla} = w^*\gamma^{\nabla}$ for some $w$,
  then $\gamma_1 = w\gamma$.

%% file: 5DynkinExt.tex
\section{\bf Extensions}

\subsection{The vertex extension}
  \label{sec_vert_ext}
     Let $\Gamma$ be one of  simply-laced Dynkin diagrams,
     and let $\Gamma'$ be one of the Carter diagram out of the homogeneous class $C(\Gamma)$.
     Each Carter diagram $\Gamma'$ from $C(\Gamma)$ has the same rank as $\rank(\Gamma)$.
     Consider any $\Gamma$-set $S = \{\tau_1, \dots, \tau_n\}$ and
     any $\Gamma'$-set, $S' = \{\tau'_1, \dots, \tau'_n\}$ obtained by
     the transformation constructed according to
     Theorem \ref{th_invol}. Then the linear spans for $S$ and $S'$ coincide:
\begin{equation}
  \label{eq_span}
     L(S) := \Sp(S), \quad L(S') := \Sp(S'), \quad L(S) = L(S').
\end{equation}
     In this chapter, we will use definition $L :=  L(S) = L(S')$.
     Let $\widetilde\Gamma$ be a Dynkin diagram with the root system $\varPhi(\widetilde\Gamma)$ such that
\begin{equation*}
   \begin{split}
      (1) & \quad \rank(\varPhi(\widetilde\Gamma)) = \rank(\varPhi(\Gamma)) + 1, \\
      (2) & \quad L \subset \Sp(\varPhi(\widetilde\Gamma)).
   \end{split}
\end{equation*}
     The choose of the root system $\widetilde\Gamma$ is ambiguous.
     The pair $\{\Gamma', \widetilde\Gamma\}$, where $\Gamma' \in C(\Gamma)$,
     is said to be the {\it vertex extension}  of $\Gamma'$ and is denoted by
     $\Gamma ' \prec \widetilde\Gamma$.

\begin{figure}[h]
\centering
\includegraphics[scale=0.24]{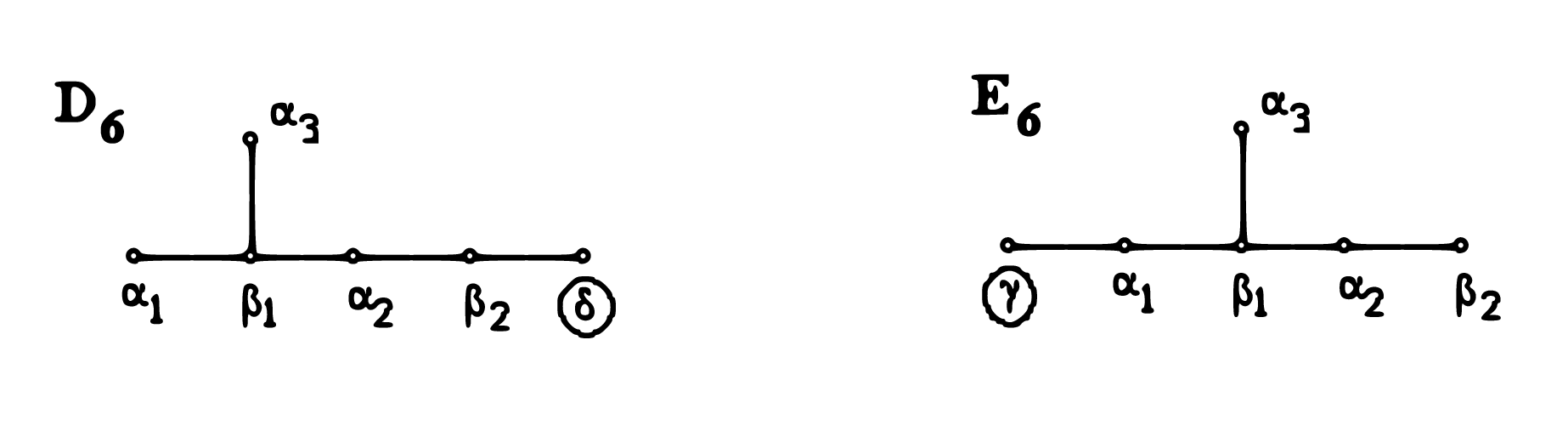}
 \caption{Two vertex extensions of $D_5$: $D_5 \prec E_6$ and $D_5 \prec D_6$.}
\label{ext_E6_D6}
\end{figure}

\begin{remark}{\rm
  It may happen that a given root subset $L$ can be extended to two different subsets
  $L(\gamma) \neq L(\delta)$
  such that $\mathscr{B}^{\vee}_{\Gamma}(\gamma^{\nabla}) \neq \mathscr{B}^{\vee}_{\Gamma}(\delta^{\nabla})$.
  For example,  there are two extensions for $D_5$ with different values of $\mathscr{B}^{\vee}_{\Gamma}(\theta^{\nabla})$:
  \begin{equation*}
    \mathscr{B}^{\vee}_{\Gamma}(\gamma^{\nabla}) = \frac{5}{4} \quad \text{ for } \quad  D_5 \prec E_6,  \qquad
    \mathscr{B}^{\vee}_{\Gamma}(\delta^{\nabla}) = 1 \quad \text{ for } \quad  D_5 \prec D_6,
  \end{equation*}
  see Fig.~\ref{ext_E6_D6}.
  The matrices $B_\Gamma$ and $B_\Gamma^{-1}$ for $\Gamma = D_5$
  are as follows:
 \begin{equation*}
  \begin{array}{c}
  B_\Gamma =
    \left [
    \begin{array}{ccccc}
    2 & 0 & 0  &   -1 &  0 \\
     0 & 2 & 0  &   -1 & -1 \\
     0 & 0 & 2  &   -1 &  0 \\
    -1 & -1 & -1 &   2 &  0 \\
     0 & -1 & 0  &   0 & 2  \\
    \end{array}
    \right ]
    \begin{array}{c}
      \alpha_1 \\
      \alpha_2 \\
      \alpha_3 \\
      \beta_1 \\
      \beta_2 \\
    \end{array}
 \end{array},
   \qquad
  B_\Gamma^{-1} =
  \begin{array}{c}
    \frac{1}{4} \left
    [ \begin{array}{ccccc}
     5 & 4 & 3 & 6 & 2  \\
     4 & 8 & 4 & 8 & 4  \\
     3 & 4 & 5 & 6 & 2  \\
     6 & 8 & 6 & 12 & 4 \\
     2 & 4 & 2 & 4 & 4 \\
    \end{array} \right ]
      \begin{array}{c}
      \alpha_1 \\
      \alpha_2 \\
      \alpha_3 \\
      \beta_1 \\
      \beta_2 \\
    \end{array}
  \end{array}
 \end{equation*}
}
\end{remark}

\subsection{Partial linkage systems}  Let $\widetilde\Gamma$ be a Dynkin diagram such that
the pair  $\{\widetilde\Gamma, \Gamma\}$ is the vertex extension $\Gamma' \prec \widetilde\Gamma$,
see \S\ref{sec_vert_ext}. Denote by $\mathscr{L}_{\widetilde\Gamma}(\Gamma', S')$ the set of
linkage diagrams
\begin{equation}
  \label{eq_part_1}
      \mathscr{L}_{\widetilde\Gamma}(\Gamma', S') =
      \{ \gamma^{\nabla} \mid \gamma \in \varPhi(\widetilde\Gamma), \gamma \not\in L(S') \},
\end{equation}
where $\varPhi(\widetilde\Gamma)$ is the root system associated with $\widetilde\Gamma$,
$S'$ some $\Gamma'$-set and
the linear space $L = L(S')$ is given by \eqref{eq_span}.
\begin{proposition}
 \label{prop_not_dep_S}
   The set $\mathscr{L}_{\widetilde\Gamma}(\Gamma', S')$ does not depend on choosing $\Gamma'$-set $S'$.
\end{proposition}
  \PerfProof
 Let $S'_1$ and $S'_2$ be two different $\Gamma'$-sets in $\varPhi(\widetilde\Gamma)$,
 and let $S_i' = \{\tau_1^i,\dots, \tau_l^i \}$ for $i = 1,2$.
 Consider vector $\gamma_1 \in \varPhi(\widetilde\Gamma)$ obtained by the vertex extension of $\Gamma'$
 using the $\Gamma'$-set $S'_1$.  Let $\gamma_1,L$ be the projection of $\gamma$ onto $L = L(S_1')$.   
 By \eqref{eq_theta_decomp} we have
\begin{equation}
      \gamma_1 =   \mu_1 + \sum\limits_{j=1}^l t_j\tau^1_j, \\
\end{equation}
where $\{t_j \mid j = 1\dots l\}$ are coefficients of decomposition of $\gamma_1,L$ by the basis 
$\{\tau_1^i,\dots, \tau_l^i \}$. 
From $\gamma_1$ we construct $\gamma_2 \in \varPhi(\widetilde\Gamma)$ by the same vertex extension using
$\Gamma'$-set $S'_2$:
\begin{equation}
       \gamma_2 =  \mu_2 + \sum\limits_{j=1}^l t_j\tau^2_j,  \\
\end{equation}
where that basis vectors $\tau_j^1$ and $\tau_j^2$ have the same coefficient $t_j$.
Here, $\mu_1, \mu_2$ are normal extending vectors from \S\ref{sec_normal_mu}.
Coordinates of the linkage label vectors $\gamma^{\nabla}_1$ and $\gamma^{\nabla}_2$ are as follows;
 \begin{equation}
     (\gamma^{\nabla}_1)_k = \sum\limits_{j=1}^l t_j(\tau^1_j, \tau^1_k), \quad
     (\gamma^{\nabla}_2)_k = \sum\limits_{j=1}^l t_j(\tau^2_j, \tau^2_k), \quad k = 1,\dots,l.
 \end{equation}
 Since inner products $(\tau^1_j, \tau^1_k)$ and $(\tau^1_j, \tau^1_k)$ correspond to the same edge
 of $\Gamma'$ (or both are equal to $0$) we get $\gamma^{\nabla}_1 = \gamma^{\nabla}_2$.
 Thus, we get
\begin{equation}
    \mathscr{L}_{\widetilde\Gamma}(\Gamma', S'_1) = \mathscr{L}_{\widetilde\Gamma}(\Gamma', S'_2).
\end{equation}
\qed

Thus, the set of linkage diagrams \eqref{eq_part_1} can be denoted by $\mathscr{L}_{\widetilde\Gamma}(\Gamma'$).
The set of linkage diagrams $\mathscr{L}_{\widetilde\Gamma}(\Gamma')$ is said to be {\it linkage system $\Gamma'$ over $\Gamma$}
or {\it partial linkage system}.

\subsection{Full linkage systems}
 The union of all partial systems by all possible vertex extensions
is said to be the {\it linkage system} or {\it full linkage system}  of $\Gamma$:
\begin{equation*}
      \mathscr{L}(\Gamma') = \bigcup_{\Gamma ' \prec \widetilde\Gamma} \mathscr{L}_{\widetilde\Gamma}(\Gamma').
\end{equation*}

\begin{proposition}
   \label{prop_sizes}
   {\rm(i)}
   For the homogeneous Carter diagrams $\Gamma$ and  $\Gamma'$
   the sizes of the full linkage systems are the same:
\begin{equation}
  \label{eq_coincide_L}
   \mid \mathscr{L}(\Gamma) \mid \;\; = \;\; \mid \mathscr{L}(\Gamma') \mid.
\end{equation}

   {\rm(ii)}
    For the full linkage systems \eqref{eq_coincide_L},
    the estimate of the number of linkage diagrams is as follows:
\begin{equation*}
     \mid \mathscr{L}(\Gamma') \mid \; = \; \mid \mathscr{L}(\Gamma) \mid
        \; \leq \; \bigcup \big (\mid \varPhi(\widetilde\Gamma_i) \mid - \mid \varPhi(\Gamma) \mid \big ),
\end{equation*}
where the union is taken by all vertex extensions $\Gamma ' \prec \widetilde\Gamma_i$.
\end{proposition}

\PerfProof
(i) This fact follows from Proposition \ref{prop_quadr_1} and \S\ref{sec_chain}.

(ii) By \eqref{eq_coincide_L}, for any $\Gamma'$ out of the homogeneous class $C(\Gamma)$
we have
\begin{equation}
   \label{eq_estimate_L}
     \mid \mathscr{L}_{\widetilde\Gamma}(\Gamma') \; = \; \mid \mathscr{L}_{\widetilde\Gamma}(\Gamma) \mid
        \;\; \leq \;\; \mid \varPhi(\widetilde\Gamma) \mid - \mid \varPhi(\Gamma) \mid.
\end{equation}
 Then, the number of linkage diagrams in the full linkage systems can be estimated as follows:
\begin{equation*}
     \mid \mathscr{L}(\Gamma') \mid \; = \; \mid \mathscr{L}(\Gamma) \mid
        \;\; \leq \;\; \bigcup \big (\mid \varPhi(\widetilde\Gamma_i) \mid - \mid \varPhi(\Gamma) \mid \big ),
\end{equation*}
where the union is taken by all vertex extensions $\Gamma ' \prec \widetilde\Gamma_i$.
\qed

%% file: 6DTypeLinkages.tex
\section{\bf The linkage systems  $\mathscr{L}(D_l)$ and $\mathscr{L}(D_l(a_k))$.}
    \label{sec_D_lg8}

 In this section we assume that $\Gamma$ is one of the Carter diagrams $D_l$ or $D_l(a_k)$.

\subsection{$D$-components $\mathscr{L}_{D_{l+1}}(D_l)$ and $\mathscr{L}_{D_{l+1}}(D_l(a_k))$}
\begin{lemma}
  \label{lem_coinc_linkage}
    Let $D_l \prec D_{l+1}$ be some vertex extension with root systems
 \begin{equation*}
      \varPhi(D_l) = \{ \tau_1, \dots, \tau_l \}, \quad
      \varPhi(D_{l+1}) = \{ \tau_1, \dots, \tau_l, \tau_{l+1} \},
 \end{equation*}
 where $\tau_{l+1}$ is a simple positive root in $\varPhi(D_{l+1})\backslash\varPhi(D_l)$.
 Let $\varphi$ be some positive root in  $\varPhi(D_{l+1})\backslash\varPhi(D_l)$,
 and $\mu_{max}$ be maximal root in $\varPhi(D_{l+1})$.

     {\rm (i)} The vector
   \begin{equation}
     \label{eq_Dl_ext1}
      \delta = \mu_{max} - \varphi + \tau_{l+1}
   \end{equation}
   is also a root in $\varPhi(D_{l+1})\backslash\varPhi(D_l)$.

    {\rm (ii)}
    The linkage label vectors $\varphi^{\nabla}$ and $-\delta^{\nabla}$ coincide.

    {\rm(iii)}
      For any $l \geq 4$, the linkage systems  $\mathscr{L}_{D_{l+1}}(D_l)$ and $\mathscr{L}_{D_{l+1}}(D_l(a_k))$,
   where $1 \leq k \leq \big [\frac{l-2}{2} \big ]$, have $2l$ elements:
\begin{equation}
   \label{eq_size_Dk}
    \mid \mathscr{L}_{D_{l+1}}(D_l) \mid \; = \; \mid \mathscr{L}_{D_{l+1}}(D_l(a_k)) \mid \; = \; 2l.
\end{equation}
\end{lemma}

\begin{figure}[H]
\centering
\includegraphics[scale=0.35]{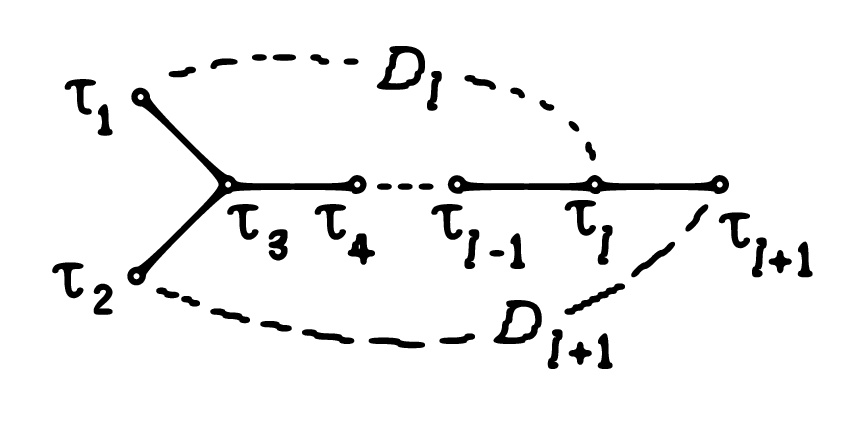}
\caption{
\hspace{3mm}The vertex extension $D_l \prec D_{l+1}$.}
\label{Dl_ext_Dl1}
\end{figure}

  \PerfProof  Let $\mathscr{B}$ be the quadratic Tits form associated with $D_{l+1}$.

  (i)  By \eqref{eq_Kac} it suffices to prove that $\mathscr{B}(\delta) = 2$.
   If $\varphi = \mu_{max}$ (resp. $\tau_{l+1}$), then $\delta = \tau_{l+1}$ (resp. $\mu_{max}$).
   In both cases, $\delta$ is the root in $\varPhi(D_{l+1})\backslash\varPhi(S)$. Suppose $\varphi \neq \mu_{max}, \tau_{l+1}$.
   We need to prove that $\mathscr{B}(\mu_{max} - \varphi + \tau_{l+1}) = 2$, i.e.,
   \begin{equation}
     \label{eq_Dl_ext3}
       \mathscr{B}(\mu_{max}) + \mathscr{B}(\varphi) + \mathscr{B}(\tau_{l+1}) +
           2(\mu_{max}, \tau_{l+1}) - 2(\varphi, \mu_{max} + \tau_{l+1}) = 2.
   \end{equation}
   Since  $\mathscr{B}(\varphi) = \mathscr{B}(\tau_{l+1}) = \mathscr{B}(\mu_{max}) = 2$,
   eq. \eqref{eq_Dl_ext3} is equivalent to the following:
  \begin{equation*}
        (\varphi, \mu_{max} + \tau_{l+1}) - (\mu_{max}, \tau_{l+1}) = 2.
   \end{equation*}
   Further, $(\mu_{max}, \tau_i) = 0$ for any $i \neq l$. In particular, $(\mu_{max}, \tau_{l+1}) = 0$ and
   it suffices to prove that
   \begin{equation}
     \label{eq_Dl_ext5}
       (\varphi, \mu_{max} + \tau_{l+1}) = 2 \text{ for any } \varphi \in \varPhi(D_{l+1})\backslash\varPhi(D_l).
   \end{equation}
   We have
   \begin{equation}
      \label{eq_Dl_ext6}
     (\gamma, \mu_{max} + \tau_{l+1}) =
       \begin{cases}
           2 \text{ for } \gamma = \tau_{l+1}, \text{ since } \mu_{max} \perp \tau_{l+1}, \\
           0 \text{ for } \gamma = \tau_l, \text{ since } (\tau_l, \mu_{max}) = 1, (\tau_l, \tau_{l+1}) = -1, \\
           0 \text{ for } \gamma = \tau_i, \text{ where } i < l.
       \end{cases}
   \end{equation}

   Since $\varphi$ is the positive root in $\varPhi(D_{l+1})\backslash\varPhi(D_l)$,
   then $\tau_{l+1}$ appears with a unit coefficient in the decomposition of $\varphi$ in the basis
   $\{\tau_1, \dots, \tau_{l+1}\}$. Then, by \eqref{eq_Dl_ext6} the relation \eqref{eq_Dl_ext5} holds for
   any root $\varphi \in \varPhi(D_{l+1})\backslash\varPhi(D_l)$.
   Therefore, $\mathscr{B}(\delta) = 2$, and $\delta$ is also a root.
   ~\\

   (ii) By \eqref{eq_Dl_ext1} we have
   \begin{equation*}
       (\delta, \tau_i) = (\mu_{max} + \tau_{l+1}, \tau_i) - (\varphi, \tau_i), \text{ where } 1 \leq i \leq l.
   \end{equation*}
   By \eqref{eq_Dl_ext6} for $i \leq l$, we have
    \begin{equation*}
     \label{ex_def_beta}
        \delta^{\nabla}_i = (\delta, \tau_i) =  -(\varphi, \tau_i) = -\varphi^{\nabla}_i, \text{ where } 1 \leq i \leq l.
   \end{equation*}
   This means that $\delta^{\nabla} = -\varphi^{\nabla}$.
   ~\\

   (iii) The number of roots of $\varPhi(D_l)$ is $2l(l-1)$, see \cite[Table IV]{Bo02}.
    Then by \eqref{eq_coincide_L} and \eqref{eq_estimate_L} we have
   \begin{equation*}
        \mid \mathscr{L}_{D_{l+1}}(D_l(a_k)) \mid ~=~ \mid \mathscr{L}_{D_{l+1}}(D_l) \mid ~\leq~ 2(l+1)l - 2l(l-1) = 4l. \\
   \end{equation*}
   By (ii), this size is two times smaller:
    \begin{equation*}
        \mid \mathscr{L}_{D_{l+1}}(D_l(a_k)) \mid ~=~ \mid \mathscr{L}_{D_{l+1}}(D_l) \mid ~\leq~ 2l.
   \end{equation*}

      In Fig.~\ref{Dlpu_linkages_2},
      we present exactly $2l$ linkage label vectors for the Carter diagram $D_l$.
      To construct the partial linkage system $\mathscr{L}_{D_{l+1}}(D_l)$  (resp. $\mathscr{L}_{D_{l+1}}(D_l(a_k))$,
      we use the action of the Weyl group $W^{\vee}_S$,
      see \S\ref{sec_partial_groups}, and Remark \ref{rem_action}, see Fig.~\ref{Dk_al_linkages}.

      So, we conclude that there are exactly $2l$ linkage label vectors in
      the $D$-component $\mathscr{L}_{D_{l+1}}(D_l)$. Then, by \eqref{eq_coincide_L} this is also true
      for the $D$-component $\mathscr{L}_{D_{l+1}}(D_l(a_k))$.  \qed

\subsection{$E$-component $\mathscr{L}_{E_8}(D_7)$.}
  By \eqref{eq_size_Dk} the linkage system $\mathscr{L}_{D_8}(D_7)$ contains $14$ elements.
  We will show that all of them lie in $\mathscr{L}_{E_8}(D_7)$.
\begin{lemma}
 \label{lem_14_pairs}
 There are $14$ pairs of roots $\eta, \lambda \in \varPhi(E_8)$ such that $\eta^{\nabla} = -\lambda^{\nabla}$.
\end{lemma}

\PerfProof
 The $7$ pairs of roots $\eta, \lambda \in \varPhi(E_8)$ such that $\eta^{\nabla} = -\lambda^{\nabla}$
listed in Table \ref{tab_pairs_roots_E8}. The opposite roots form $7$ pairs of negative roots
$\eta, \lambda \in \varPhi(E_8)$ such that $\eta^{\nabla} = -\lambda^{\nabla}$.
For $1 \leq i \leq 7$, the sum $\eta_i + \lambda_i$ is the same vector, the order of coordinates is shown
in parentheses:
\begin{equation}
   \label{eq_v_sum}
     \eta_i + \lambda_i \;=\;
   \begin{array}{ccccccc}
     4 & 7 & 10 & 8 & 6 & 4 & 2 \\
            &   & 5 \\
      \end{array}, \quad
   \biggl(
   \begin{array}{ccccccc}
     \tau_1 & \tau_3 & \tau_4 & \tau_5 & \tau_6 & \tau_7 & \tau_8 \\
            &   & \tau_2 \\
      \end{array}  \biggr)
\end{equation}
Vectors $\eta_i + \lambda_i$ from \eqref{eq_v_sum} are orthogonal to any $\tau_i, \; 2 \leq i \leq 8$, so
$\eta_i^{\nabla} = -\lambda_i^{\nabla}$. Note that there is no need orthogonality for $\tau_1$,
see \eqref{eq_dual_gamma}.  Further, we calculate the label vectors $\eta_i^{\nabla}$
(or $\lambda_i^{\nabla}$). They are listed in the last column of Table \ref{tab_pairs_roots_E8}.
 \begin{table}[!ht]
  \centering
  \renewcommand{\arraystretch}{1.3}
  \begin{tabular} {|c|c|c|c|}
  \hline
       & \footnotesize $\eta \in \varPhi(E_8)\backslash\varPhi(D_7)$ &
          \footnotesize$\lambda \in \varPhi(E_8)\backslash\varPhi(D_7)$
       & \footnotesize$\eta^{\nabla} = -\lambda^{\nabla}\in \mathscr{L}_{E_8}(D_7)$ \\
    \hline
     \footnotesize $1$ & \footnotesize$\begin{array}{ccccccc}
            2 & 3 & 4 & 3 & 2 & 1 & 0 \\
              &  &  2
         \end{array}$ &
         \footnotesize$\begin{array}{ccccccc}
            2 & 4 & 6 & 5 & 4 & 3 & 2 \\
              &  & 3
         \end{array}$ &
          \footnotesize$\begin{array}{ccccccc}
             0 & 0 & 0 & 0 & 0 & -1 \\
              & 0
         \end{array}$ \\
    \hline
     \footnotesize$2$ & \footnotesize$\begin{array}{ccccccc}
            2 & 3 & 4 & 3 & 2 & 1 & 1 \\
              & & 2
         \end{array}$ &
         \footnotesize$\begin{array}{ccccccc}
            2 & 4 & 6 & 5 & 4 & 3 & 1 \\
              &  & 3
         \end{array}$ &
         \footnotesize $\begin{array}{ccccccc}
             0 & 0 & 0 & 0 & -1 & 1 \\
              & 0
         \end{array}$ \\
    \hline
     \footnotesize$3$ & \footnotesize$\begin{array}{ccccccc}
            2 & 3 & 4 & 3 & 2 & 2 & 1 \\
              & & 2
         \end{array}$ &
         \footnotesize$\begin{array}{ccccccc}
            2 & 4 & 6 & 5 & 4 & 2 & 1 \\
              & & 3
         \end{array}$ &
          \footnotesize$\begin{array}{ccccccc}
             0 & 0 & 0 & -1 & 1 & 0 \\
              & 0
         \end{array}$ \\
    \hline
    \footnotesize$4$ & $\begin{array}{ccccccc}
            2 & 3 & 4 & 3 & 3 & 2 & 1 \\
              & & 2
         \end{array}$ &
    \footnotesize $\begin{array}{ccccccc}
            2 & 4 & 6 & 5 & 3 & 2 & 1 \\
              & & 3
         \end{array}$ &
    \footnotesize $\begin{array}{ccccccc}
             0 & 0 & -1 & 1 & 0 & 0 \\
              & 0
         \end{array}$ \\
    \hline
     \footnotesize$5$ & \footnotesize$\begin{array}{ccccccc}
            2 & 3 & 4 & 4 & 3 & 2 & 1 \\
              & & 2
         \end{array}$ &
     \footnotesize $\begin{array}{ccccccc}
            2 & 4 & 6 & 4 & 3 & 2 & 1 \\
              & & 3
         \end{array}$ &
     \footnotesize $\begin{array}{ccccccc}
             0 & -1 & 1 & 0 & 0 & 0 \\
              & 0
         \end{array}$ \\
    \hline
     \footnotesize $6$ & \footnotesize $\begin{array}{ccccccc}
            2 & 3 & 5 & 4 & 3 & 2 & 1 \\
              & & 2
         \end{array}$ &
     \footnotesize $\begin{array}{ccccccc}
            2 & 4 & 5 & 4 & 3 & 2 & 1 \\
              & & 3
         \end{array}$ &
     \footnotesize $\begin{array}{ccccccc}
             -1 & 1 & 0 & 0 & 0 & 0 \\
              & -1
         \end{array}$ \\
    \hline
     \footnotesize $7$ & \footnotesize $\begin{array}{ccccccc}
            2 & 3 & 5 & 4 & 3 & 2 & 1 \\
              & & 3
         \end{array}$ &
     \footnotesize $\begin{array}{ccccccc}
            2 & 4 & 5 & 4 & 3 & 2 & 1 \\
              & & 2
         \end{array}$ &
     \footnotesize $\begin{array}{ccccccc}
             -1 & 0 & 0 & 0 & 0 & 0 \\
                & 1
         \end{array}$ \\
    \hline
       \end{tabular}
  \vspace{2mm}
  \caption{The $7$ pairs of positive roots $\eta, \lambda \in \varPhi(E_8)$
  such that $\eta^{\nabla} = -\lambda^{\nabla}$.}
  \label{tab_pairs_roots_E8}
  \end{table}

\subsection{The size and structure of linkage systems $\mathscr{L}(D_l)$ and  $\mathscr{L}(D_l(a_k))$}
  \label{sec_size_n_str}

\begin{theorem}
     \label{th_size_Dtype}

    {\rm(i)}
      For $l = 4$, the full linkage  system $\mathscr{L}(D_4)$ (resp. $\mathscr{L}(D_4(a_1))$
      consists only of $D$-component, see Figs.~\ref{D4_loctets},~\ref{D4a1_loctets}.
      Each of these $D$-components consists of three orbits of size $8$ each.

     {\rm(ii)}
      For $l > 7$, the full linkage system $\mathscr{L}(D_l)$ (resp. $\mathscr{L}(D_l(a_1))$)
      consists only of $D$-component for any $k$, see Figs.~\ref{Dlpu_linkages_2},~\ref{Dk_al_linkages}.
      Each of these $D$-components consists of one orbit of size $2l$.

    {\rm(iii)} For $l=5$, the full linkage system $\mathscr{L}(D_5)$
    (resp. $\mathscr{L}(D_5(a_1))$) consists of the $D$-component ($10$ elements) and
    the $E$-component consisting of two orbits ($2 \times 16$ elements).
    In total, $\mathscr{L}(D_5)$  contains $42$ elements.
    One of two orbits of the $E$-component of the $\mathscr{L}(D_5)$
    is shown in Fig.~\ref{D5pu_loctets}.

    {\rm(iv)} For $l=6$, the full linkage system $\mathscr{L}(D_6)$
    (resp. $\mathscr{L}(D_6(a_1))$, resp. $\mathscr{L}(D_6(a_2))$ consists of the $D$-component ($12$ elements) and
    the $E$-component consisting of two orbits ($2 \times 32$ elements).
    In total, $\mathscr{L}(D_6)$ contains $76$ elements.
    Two orbits of the $E$-component of the linkage system $\mathscr{L}(D_6)$
    are shown in Fig.~\ref{D6pu_loctets}.

    {\rm(v)} For $l=7$, the full linkage system $\mathscr{L}(D_7)$ (resp. $\mathscr{L}(D_7(a_1))$,
    resp. $\mathscr{L}(D_7(a_2))$) consists of one $D$-component ($14$ elements) and
    two $E$-components consisting of two orbits ($2 \times 64$ elements).
    In total, $\mathscr{L}(D_7)$ contains $142$ elements.
  \end{theorem}

\PerfProof
  (i) Three orbits of $\mathscr{L}_{D_5}(D_4)$ are shown in Fig.~\ref{D4_loctets}.
   Three orbits $\mathscr{L}_{D_5}(D_4(a_1))$ are shown in Fig.~\ref{D4a1_loctets}.

  (ii) For $l > 7$, the $D$-component $\mathscr{L}_{D_{l+1}}(D_l)$ and $\mathscr{L}_{D_{l+1}}(D_l(a_k))$
   are depicted in Figs.~\ref{Dlpu_linkages_2} and \ref{Dk_al_linkages}.
   By Lemma \ref{lem_coinc_linkage}(iii) there are $2l$ elements in each of these components.

  (iii) By \eqref{eq_coincide_L} and \eqref{eq_estimate_L} for $l= 5$,
   \begin{equation*}
     \mid \mathscr{L}_{E_6}(D_5(a_1)) \mid ~=~ \mid \mathscr{L}_{E_6}(D_6) \mid
      ~\leq~ \mid \varPhi(E_6) \mid -  \mid \varPhi(D_5) \mid ~=~ 72 - 40 = 32. \\
 \end{equation*}
   The full linkage system $\mathscr{L}(D_5)$ consists of $10$-element linkage system $\mathscr{L}_{D_6}(D_5)$
   (Lemma \ref{lem_coinc_linkage}(iii)) and two orbits of the $E$-component, each of which contains
   $16$ linkages, i.e.,
 \begin{equation}
   \label{eq_l_D5}
     \mid \mathscr{L}(D_5) \mid ~=~ \mid \mathscr{L}(D_5(a_1)) \mid ~\leq~ 32 + 10 = 42. \\
 \end{equation}
   One of these components is depicted in Fig.~\ref{D5pu_loctets}. Another $E$-component
   is constructed from the first one as follows: for any linkage diagram the solid edges become dotted,
   and vice versa, the corresponding label vector becomes just opposite.
   Thus, in \eqref{eq_l_D5} both
   $\mathscr{L}(D_5)$ and $\mathscr{L}(D_5(a_1))$ contain exactly $42$ elements.
~\\

  (iv) Similarly, for $l= 6$,
   \begin{equation*}
     \mid \mathscr{L}_{E_7}(D_6(a_2)) \mid ~=~ \mid \mathscr{L}_{E_7}(D_6(a_1)) \mid ~=~
     \mid \mathscr{L}_{E_7}(D_6) \mid ~\leq~ \mid \varPhi(E_7) \mid -  \mid \varPhi(D_6) \mid ~\leq~ 126 - 60 = 66. \\
 \end{equation*}
 For the full linkage system, we have the following estimate:
   \begin{equation*}
     \mid \mathscr{L}(D_6) \mid ~=~ \mid \mathscr{L}(D_6(a_2)) \mid ~=~ \mid \mathscr{L}(D_6(a_1)) \mid
      ~\leq~ 66 + 12 = 78. \\
 \end{equation*}
 Let the coordinates of roots of $E_7$ be as follows
  \begin{equation*}
   \begin{array}{cccccc}
     \tau_6 & \tau_1 & \tau_2 & \tau_3 & \tau_4 & \tau_5 \\
            &   & \beta_2 \\
    \end{array}
   \end{equation*}
 Consider maximal and minimal roots $\pm\mu$ in $E_7$. Here,
 $\pm\mu \in \varPhi(E_7) \backslash \varPhi(D_6)$ and $\pm\mu_{max}\not\in \varPhi(D_6)$.
 Note that $\pm\mu$ are orthogonal to any simple root except for $\tau_6$.
 Then vectors $\pm\mu^{\nabla}$ are zero, i.e.,
   \begin{equation}
    \label{eq_l_D6}
     \mid \mathscr{L}(D_6) \mid ~=~ \mid \mathscr{L}(D_6(a_1)) \mid ~=~ \mid \mathscr{L}(D_6(a_2)) \mid
      ~\leq~ 78 - 2 = 76. \\
   \end{equation}

  The full linkage system for $\mathscr{L}(D_6)$ consists of the $12$-element linkage system $\mathscr{L}_{D_7}(D_6)$
  (Lemma \ref{lem_coinc_linkage}(iii)) and two orbits of the $E$-component, each of which contains $32$ linkages.
  These $E$-components are depicted in Fig.~\ref{D6pu_loctets},
  they contain $2 \times 32 = 64$ elements. Thus, $D$-component and $E$-components together
  contain $12 + 64 = 76$ elements in $\mathscr{L}(D_6)$.
  So, in \eqref{eq_l_D6} all full linkage systems $\mathscr{L}(D_6)$, $\mathscr{L}(D_6(a_1))$ and $\mathscr{L}(D_6(a_1))$
  contain exactly $76$ elements.
~\\

 (v) For $l= 7$, we have
  \begin{equation*}
     \mid \mathscr{L}_{E_8}(D_7(a_2)) \mid ~=~ \mid \mathscr{L}_{E_8}(D_7(a_1)) \mid ~=~
     \mid \mathscr{L}_{E_8}(D_7) \mid ~\leq~ \mid \varPhi(E_8) \mid -  \mid \varPhi(D_7) \mid
     ~\leq~ 240 - 84 = 156. \\
 \end{equation*}
 By Lemma \ref{lem_14_pairs}(i), in $\mathscr{L}_{E_8}(D_7)$ there are $14$ pairs having the same
 linkage label vectors, then

  \begin{equation*}
     \mid \mathscr{L}_{E_8}(D_7(a_2)) \mid ~=~ \mid \mathscr{L}_{E_8}(D_7(a_1)) \mid ~=~
     \mid \mathscr{L}_{E_8}(D_7) \mid ~\leq~ 156 - 14 = 142. \\
 \end{equation*}
\smallskip

In fact, $\mid \mathscr{L}_{E_8}(D_7) \mid$ contains only $128$ elements in two $64$-element orbits
which can be found in \cite[Figs. C.63, C.64]{St14}.
  The full linkage system for $\mathscr{L}(D_7)$ consists of the $14$-element linkage system $\mathscr{L}_{D_8}(D_7)$
  (Lemma \ref{lem_coinc_linkage}(iii)) and two orbits of the $E$-component.
  The full linkage systems $\mathscr{L}(D_7)$, $\mathscr{L}(D_7(a_1))$
  and $\mathscr{L}(D_7(a_2))$ contain exactly $142$ elements. \qed

%% file: 7ExamplesDiagr.tex

\begin{figure}[!ht]
\centering
\includegraphics[scale=0.24]{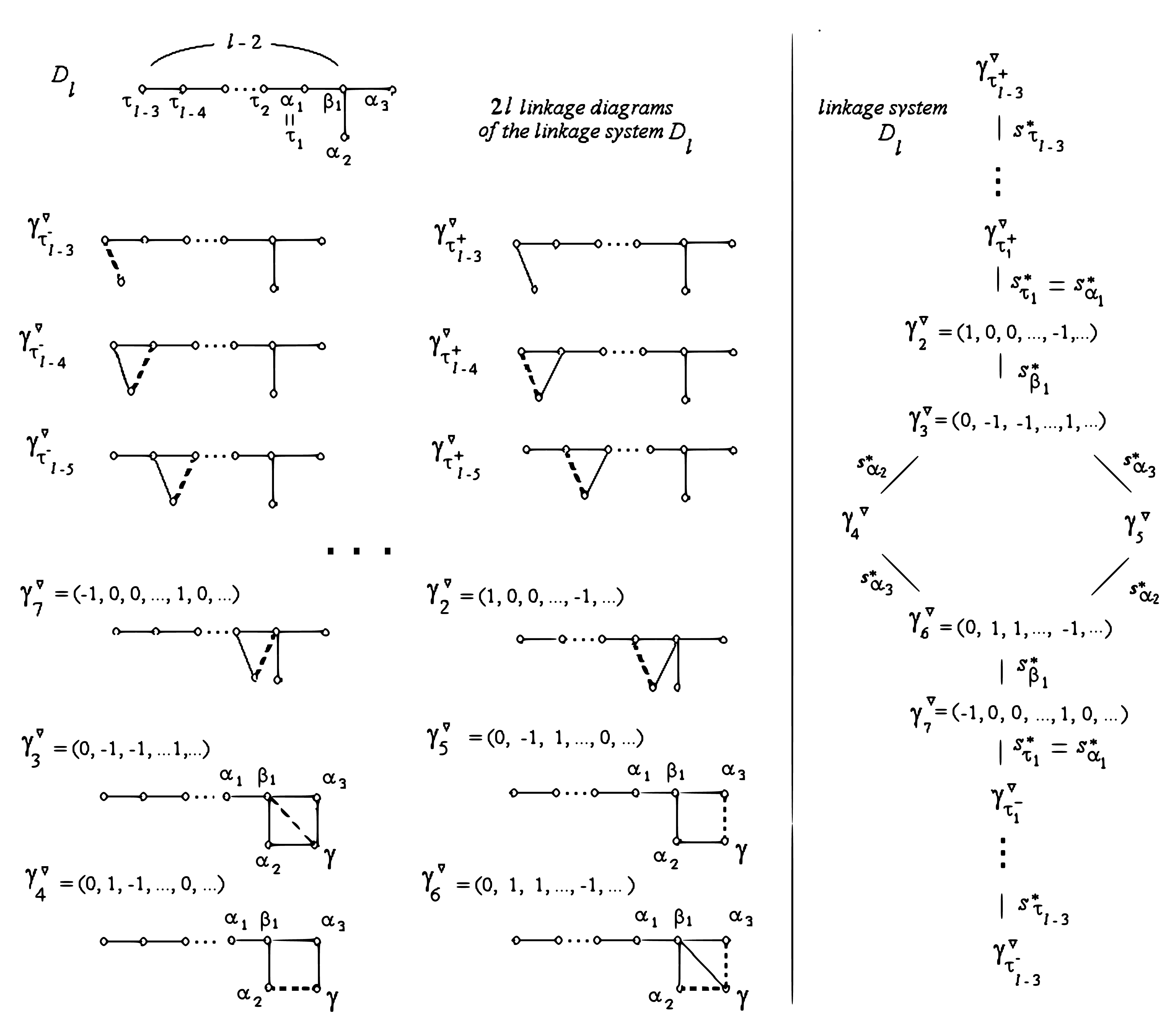}
\caption{
For $l > 4$, the partial linkage system $\mathscr{L}_{D_{l+1}}(D_l)$ contains $2l$ linkages.
For $l > 7$, $\mathscr{L}_{D_{l+1}}(D_l)$ coincides with the full linkage system $\mathscr{L}(D_l)$.}
\label{Dlpu_linkages_2}
\end{figure}

\begin{figure}[!ht]
\centering
\includegraphics[scale=0.20, angle=90]{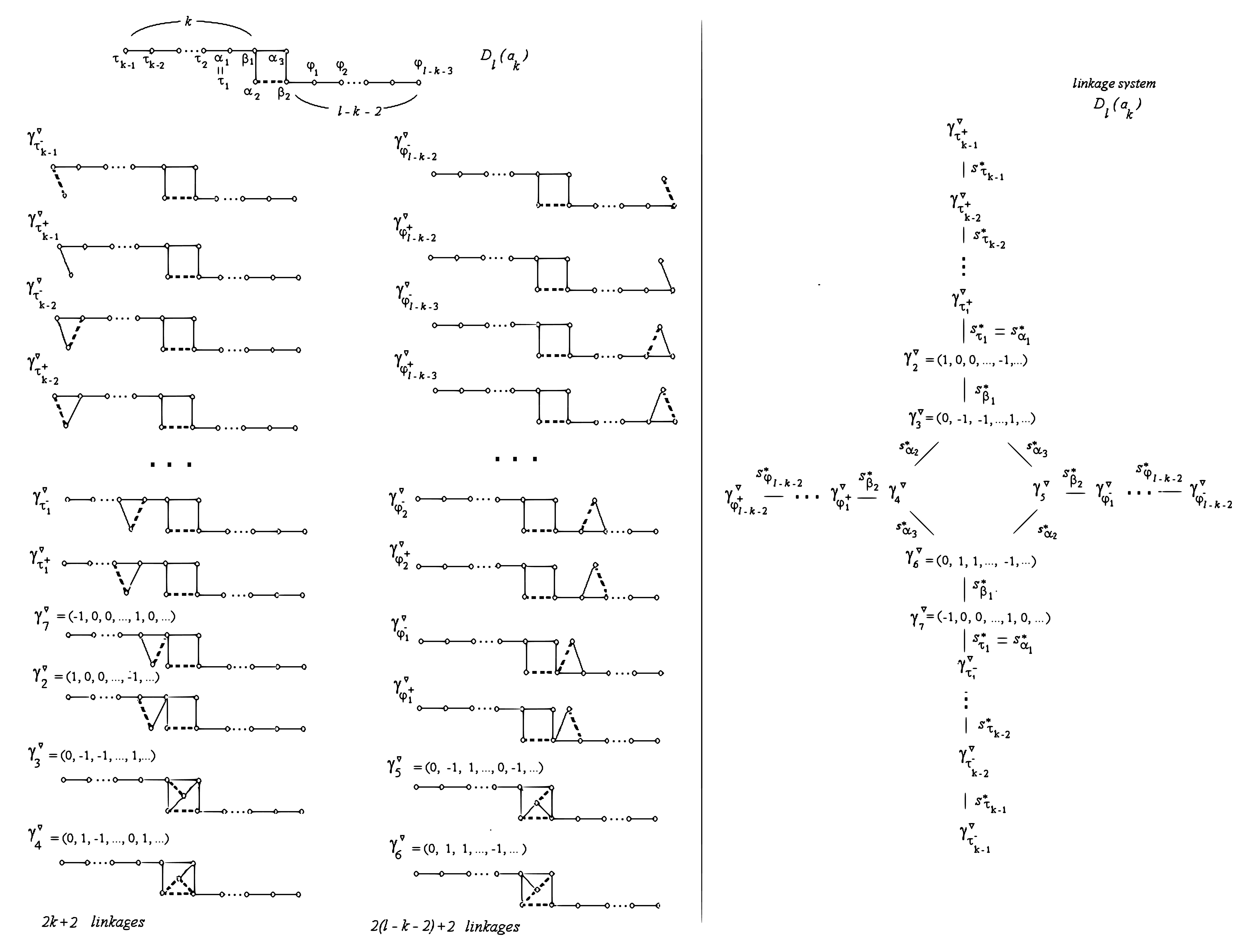}
\caption{For $l > 4$, the partial linkage system $\mathscr{L}_{D_{l+1}}(D_l(a_k))$ contains $2l$ linkages.
For $l > 7$, $\mathscr{L}(D_l(a_k)) = \mathscr{L}_{D_{l+1}}(D_l(a_k))$ .}
\label{Dk_al_linkages}
\end{figure}

\begin{figure}[!ht]
\centering
\includegraphics[scale=0.25]{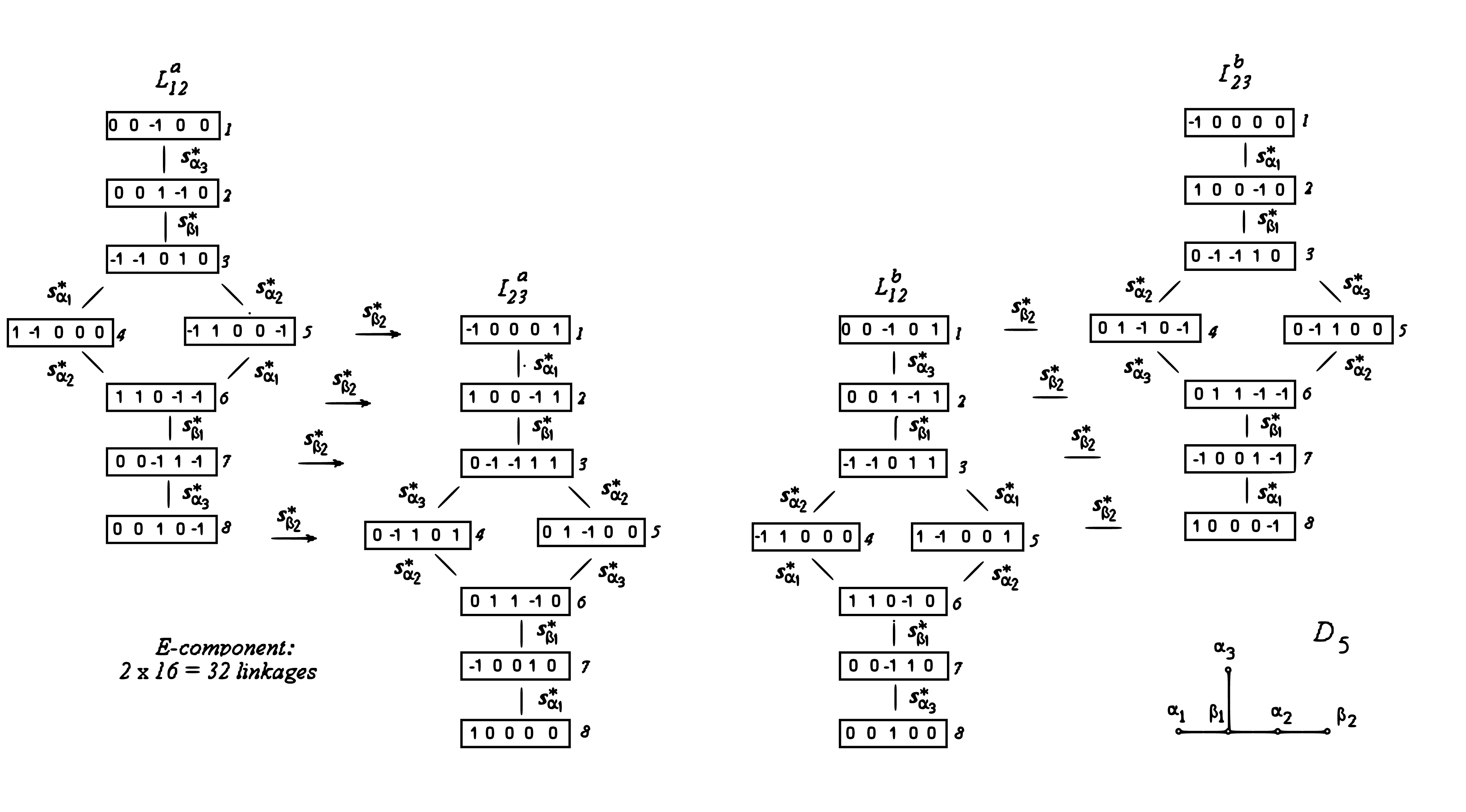}
\caption{Two orbits of the $E$-component of the linkage system $\mathscr{L}(D_5)$. 
Each orbit consists of two loctets.}
\label{D5pu_loctets}
\end{figure}

\begin{figure}[!ht]
\centering
\includegraphics[scale=0.2]{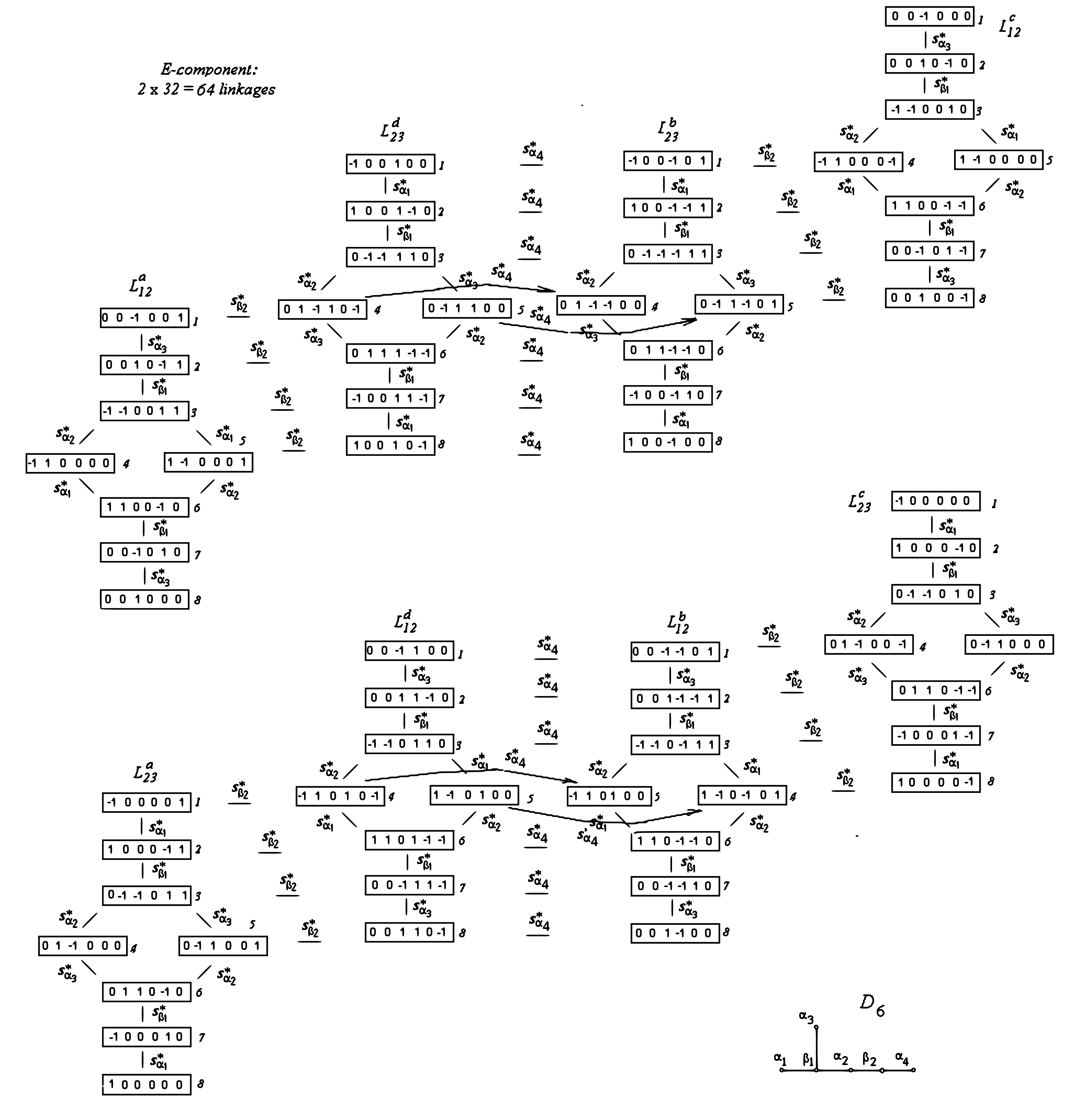}
\caption{Two orbits of the $E$-component of the linkage system $\mathscr{L}(D_6)$.
Each orbit consists of four loctets.}
\label{D6pu_loctets}
\end{figure}

%% file: App_Gabrielov.tex
\clearpage
\section{\bf Gabrielov transformations and the Ovsienko theorem}
 \subsection{Gabrielov transformations}
 \label{sec_Gabrielov}
\begin{figure}[!ht]
\centering
\includegraphics[scale=0.3]{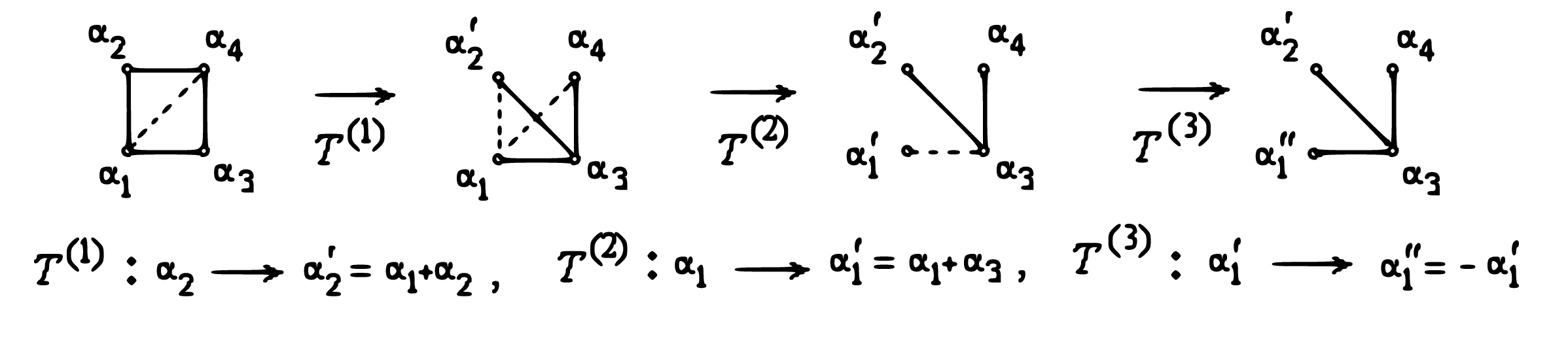}
 \caption{Gabrielov's example of changing the basis corresponding to a singularity
 $x^3 + y^3 + z^2$.
 Transitions $T^{(i)}, i = 1,2,3$ act only on one basis vector and
 fix all others.}
\label{fig_D4a1_Gabrielov}
\end{figure}

In 1973, A.~Gabrielov introduced mappings changing quadratic forms associated with
singularities, \cite[\S6]{G73}. Gabrielov's example corresponding to the singularity $x^3 + y^3 + z^2$ is shown in Fig.~\ref{fig_D4a1_Gabrielov}.
Essentially, Gabrielov transformations in Fig.~\ref{fig_D4a1_Gabrielov}
similar to transitions of \S\ref{sec_trans_mappings}.
The sequence of transformation $T^{(i)}, i = 1,2,3$  modify the underlying diagram to the
Dynkin diagram $D_4$, change the basis (vanishing cycles) and the matrix of the quadratic form
(intersection matrix) which is similar to a certain Cartan matrix, [E18], [Go99].
Gabrielov transformations are used for classification of some positive definite forms in \cite{DD95}.

\begin{figure}[!ht]
\centering
\includegraphics[scale=0.3]{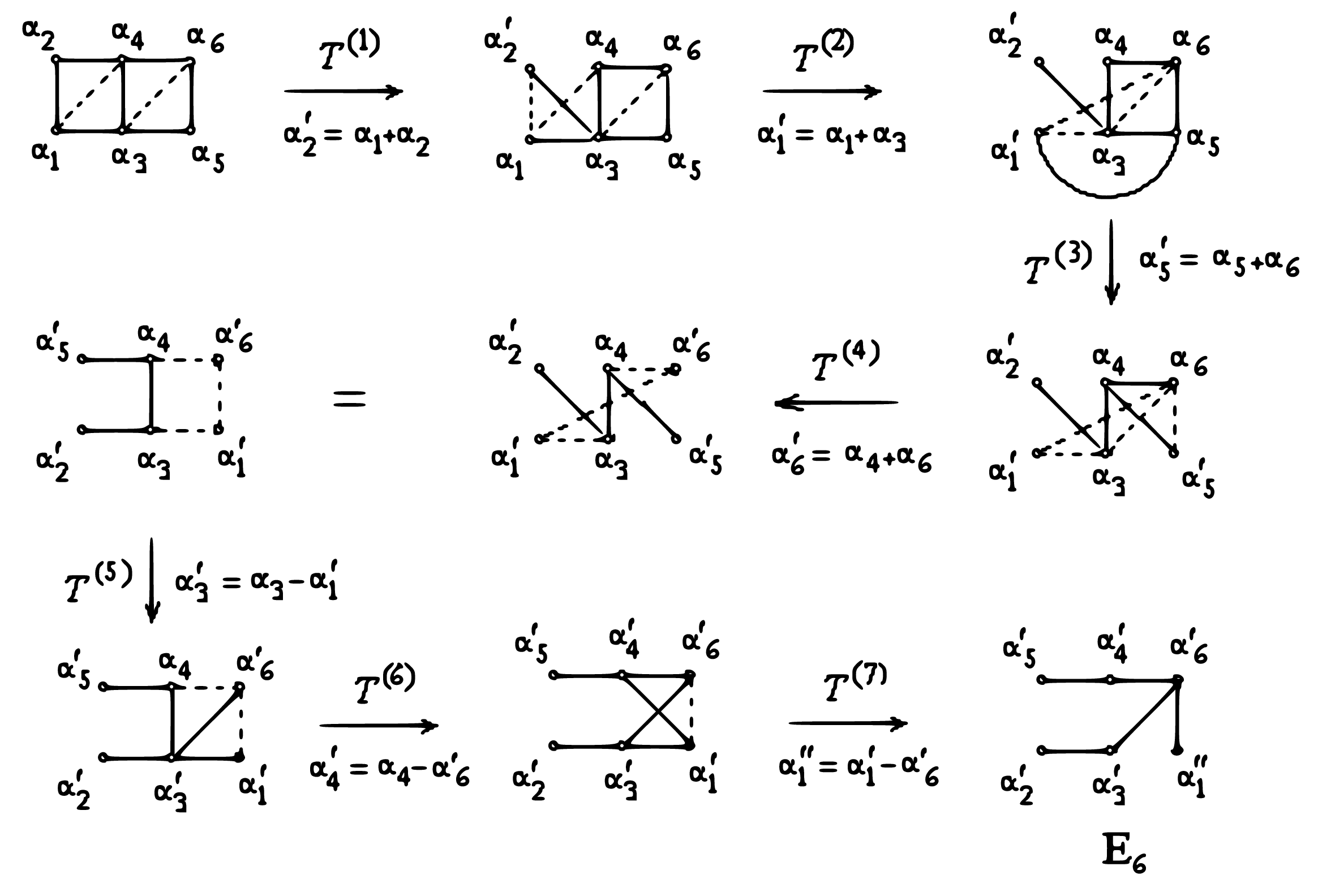}
 \caption{ The Gabrielov example of changing the basis corresponding to a singularity $x^4 + y^3 + z^2$.
 Transitions $T^{(i)}, 1 \leq i \leq 7$ act only on one basis vector and
 fix all others.}
\label{fig_E6a1_Gabrielov}
\end{figure}
Another example corresponds to the singularity $x^4 + y^3 + z^2$, see Fig.~\ref{fig_E6a1_Gabrielov}.
The sequence of Gabrielov transformations $T^{(i)}, 1 \leq i \leq 7$ in Fig. \ref{fig_E6a1_Gabrielov}
modify the underlying diagram to the Dynkin diagram $E_6$.

 Let $B = B_\Gamma$ be the (partial) Cartan matrix associated with the diagram $\Gamma$,
 $E_{ij}$ be the elementary $n \times n$ matrix with unique nonzero entry a $1$ at slot $(i,j)$,
 and $G^B_{ij}$ - the linear transformation given by matrix $G^B_{ij} = I - b_{ij}E_{ij}$.
 For each $i \neq j$ in $\{1,\dots,n \}$, the Gabrielov transformation $\mathscr{G}_{ij}$ is defined as
 the following composition:
\begin{equation}
  \label{eq_Gabrielov_tr}
   B \longrightarrow \mathscr{G}_{ij}(B) = B G^B_{ij},
\end{equation}
see \cite[Prop.~2.17]{BGP19}.

\subsection{Inflations and Ovsienko's theorem}
 \label{sec_Ovsienko}
 Let $\varepsilon$ be a sign  $\varepsilon \in \{+, - \}$, and for $1 \leq i,j \leq  n$ define the
 linear transformation $T^{\varepsilon_{ij}} :\mathbb{Z}^n \longmapsto \mathbb{Z}^n$ by
\begin{equation*}
     T^{\varepsilon}_{ij} : v  \longmapsto  v - \varepsilon v_i \alpha_j.
\end{equation*}
Note that $T^{+}_{ij}$ is the inverse of $T^{-}_{ij}$.
The transformation $T^{-}_{ij}$ is called  a {\it deflation} for $B_{\Gamma'}$
if $b_{ij} < 0$, and   $T^{+}_{ij}$ is called an {\it inflation} for $B_{\Gamma'}$ if $b_{ij} > 0$.
Inflations and deflations are called {\it flations}.
A finite composition of {\it flations} is an {\it iterated flation}:
\begin{equation}
  \label{eq_iterated_fl}
  T = T^{\varepsilon_1}_{i_1 j_1}\dots T^{\varepsilon_r}_{i_r j_r}.
\end{equation}
Composition \eqref{eq_iterated_fl} is inductively defined, see \cite[\S2.4]{BGP19}.

The inflations technique introduced by
S.~Ovsienko in \cite{O78} in the context of weakly positive unit quadratic forms,
see \cite[Rem.~4.4]{K05} and \cite[Exam.~5.2]{MM19}, is similar to Gabrielov's transformations.
If $\abs{b_{ij}} = 1$  and $\varepsilon \in \{+, - \}$ is such that $\abs{b_{ij}} = \varepsilon b_{ij}$
then the quadratic form $B T^{\varepsilon}_{ij}$ coincides with the Gabrielov transformation
$\mathscr{G}_{ij}(B)$ from \eqref{eq_Gabrielov_tr}.

\begin{theorem-non}[Ovsienko]
  Let $B$ be a positive definite unit form. Then there exists an iterated inflation
  $T$ and a unique (up to permutation of components) disjoint union of Dynkin diagrams $\Gamma$ such that
  $BT$ coincides with the Cartan matrix $B_{\Gamma}$.
\end{theorem-non}
~\\
For the proof, see  \cite[Th.~2.20]{BGP19}.
For the further development of the inflation mappings, see
\cite{BP99}, \cite{H98}, \cite{K05}. \cite{MM19}, \cite{MSZ17}.